\documentclass[a4paper, 11pt]{amsart}

\usepackage{graphics, graphicx, multicol}
\usepackage{color}
\usepackage{amssymb,amsbsy,amsmath,amsfonts,amssymb,amscd}
\usepackage{latexsym,euscript,epsfig}
\usepackage{pstricks}
\usepackage{array,delarray}
\usepackage{pst-node}
\usepackage{pspicture}
\usepackage[matrix,ps,xdvi]{xypic}

\usepackage{algorithm}
\usepackage{algorithmic}

\theoremstyle{plain}
   \newtheorem{theorem}{Theorem}[section]
   \newtheorem{proposition}[theorem]{Proposition}
   \newtheorem{lemma}[theorem]{Lemma}

\theoremstyle{definition}
   \newtheorem{definition}[theorem]{Definition}
   \newtheorem{example}[theorem]{Example}

   \newtheorem{remark}[theorem]{Remark}
   \newtheorem{note}[theorem]{Note}

\numberwithin{equation}{section}

\newcommand{\fominPFibo}{\widehat{P}}
\newcommand{\fominQFibo}{\widehat{Q}}

\newcommand\linearext[1]{Ext(#1)}

\def\permOneSymb{{\bf x}}
\def\youngLattice{\mathbb{YL}}
\def\youngFiboLattice{\mathbb{YFL}}
\def\youngFiboTableau{\mathbb{YFT}}
\def\standardYoungFiboTableau{\mathbb{SYFT}}
\newcommand{\fiboTableauNumber}[1]{\mathcal{F}_{#1}}
\newcommand{\fiboClass}[1]{\mathbb{YFC}(#1)}
\newcommand{\invNumber}[1]{\#inv(#1)}
\newcommand{\invSet}[1]{inv(#1)}
\newcommand{\ordSet}[1]{ord(#1)}
\newcommand{\evacuatedTableau}[1]{ev(#1)}
\newcommand{\evacuatedTableauFromLetter}[2]{ev(#1,#2)}

\newcommand\youngFiboTab[1]{\mathbb{YFT}_{#1}}
\newcommand\youngTab[1]{\mathbb{YT}_{#1}}
\def\younFiboTableauOrder{ \preceq}
\def\younTableauWeakOrder{ \preceq_{weak}}
\def\younTableauChainOrder{ \preceq_{chain}}
\def\partitDominanceOrder{ \geq_{dom}}
\def\permutohedronOrder{\leq_{p}}
\newcommand\minCano[1]{w_{min}^{#1}}
\newcommand\maxCano[1]{w_{max}^{#1}}
\newcommand\transposition[1]{\delta_{#1}}
\newcommand\canoposet[1]{\mathbb{P}_{#1}}

\def\posetp{\mathbb{P}}

\def\covers{\rhd}

\def\GrTeXBox#1{\vbox{\vskip\vcadre\hbox{\hskip\hcadre%
      $#1$%
   \hskip\hcadre}\vskip\vcadre}}
\def\arx#1[#2]{\ifcase#1 \relax \or%
  \ar @{-}[#2]  \or%
  \ar @2{-}[#2] \or%
  \ar @{--}[#2] \or%
  \ar @2{.}[#2] \or%
  \ar @{~}[#2]  \fi}

    \def\Affect{\, := \,}

\def\snakeTwoTwoOneTwo{
\begin{psmatrix}[colsep=\colSep,rowsep=\rowSep] \corner &
\corner &\ & \square\\ \square & \square & \corner & \square
\end{psmatrix} }

\def\snakeTwoTwoTwo{
\begin{psmatrix}[colsep=\colSep,rowsep=\rowSep] \corner &
\corner & \corner\\ \square & \square & \square
\end{psmatrix} }

\def\snakeTwoOneOneTwo{
\begin{psmatrix}[colsep=\colSep,rowsep=\rowSep] \corner &
\ &\ & \square\\ \square & \corner & \square & \square
\end{psmatrix} }

\def\snakeTwoOneTwo{
\begin{psmatrix}[colsep=\colSep,rowsep=\rowSep] \corner &\ & \square\\
\square & \corner & \square \end{psmatrix} }

\def\snakeOneOneOneOne{
\begin{psmatrix}[colsep=\colSep,rowsep=\rowSep] \corner & \square & \square &
\square \end{psmatrix} }

\def\snakeTwoOneOne{
\begin{psmatrix}[colsep=\colSep,rowsep=\rowSep] \corner\\ \square & \corner & \square
\end{psmatrix} }

\def\snakeOneTwoOne{
\begin{psmatrix}[colsep=\colSep,rowsep=\rowSep] & \square\\ \corner & \square & \square
\end{psmatrix} }

\def\snakeTwoTwo{
\begin{psmatrix}[colsep=\colSep,rowsep=\rowSep] \corner &
\corner\\ \square & \square \end{psmatrix} }

\def\snakeOneOneTwo{
\begin{psmatrix}[colsep=\colSep,rowsep=\rowSep] & & \square\\
\corner & \square & \square \end{psmatrix} }

\def\snakeOneOneOne{
\begin{psmatrix}[colsep=\colSep,rowsep=\rowSep] \corner & \square & \square \end{psmatrix} }

\def\snakeTwoOne{
\begin{psmatrix}[colsep=\colSep,rowsep=\rowSep] \corner\\ \square &
\corner \end{psmatrix} }

\def\snakeOneTwo{
\begin{psmatrix}[colsep=\colSep,rowsep=\rowSep] & \square\\ \corner &
\square \end{psmatrix} }

\def\snakeOneOne{
\begin{psmatrix}[colsep=\colSep,rowsep=\rowSep] \corner &
\square\\ \end{psmatrix} }

\def\hspace{}

\def\snakeTwo{
\begin{psmatrix}[colsep=\colSep,rowsep=\rowSep] \corner\\
\square \hspace \end{psmatrix} }

\def\snakeOne{
\begin{psmatrix}[colsep=\colSep,rowsep=\rowSep] \corner\end{psmatrix} }

\def\emptysnake{
\begin{psmatrix}[colsep=\colSep,rowsep=\rowSep] \emptyset\end{psmatrix} }

\def\yfClassUn{{\small
\begin{psmatrix}[colsep=\matrixColSep,rowsep=0] \fbox{1234} & ^{\rho = 0}  \end{psmatrix}} }

\def\yfClassCinq{{\small
\begin{psmatrix}[colsep=\matrixColSep,rowsep=0] \fbox{2134} \end{psmatrix}}  }

\def\yfClassTrois{{\small
\begin{psmatrix}[colsep=\matrixColSep,rowsep=0] \fbox{1324}\\3124\end{psmatrix}} }

\def\yfClassDeux{{\small
\begin{psmatrix}[colsep=\matrixColSep,rowsep=0] \fbox{1243} & ^{\rho = 1}\\1423\\4123\end{psmatrix}}   }

\def\yfClassSept{{\small
\begin{psmatrix}[colsep=\matrixColSep,rowsep=0] 2314\\ \fbox{3214}\end{psmatrix}} }

\def\yfClassSix{{\small
\begin{psmatrix}[colsep=\matrixColSep,rowsep=0] \fbox{2143}\\2413\\4213 \end{psmatrix}}  }

\def\yfClassQuatre{{\small
\begin{psmatrix}[colsep=\matrixColSep,rowsep=0] 1342 & ^{\rho = 2}\\\fbox{1432}\\4132\end{psmatrix}} }

\def\yfClassHuit{{\small
\begin{psmatrix}[colsep=\matrixColSep,rowsep=0] 2341\\ 2431\\ \fbox{4231}\end{psmatrix}} }

\def\yfClassNeuf{{\small
\begin{psmatrix}[colsep=\matrixColSep,rowsep=0] 3142 & ^{\rho = 3}\\\fbox{3412}\\4312\end{psmatrix}} }

\def\yfClassDix{{\small
\begin{psmatrix}[colsep=\matrixColSep,rowsep=0] 3241 & ^{\rho = 4}\\ 3421\\ \fbox{4321}\end{psmatrix}}  }

\def\tUn{{\small
\begin{psmatrix}[colsep=\matrixColSep,rowsep=\matrixRowSep] 5&4&3&2&1 \end{psmatrix}} }

\def\tOnze{{\small
\begin{psmatrix}[colsep=\matrixColSep,rowsep=\matrixRowSep] &&&2\\5&4&3&1 \end{psmatrix}} }
\def\tCinq{{\small
\begin{psmatrix}[colsep=\matrixColSep,rowsep=\matrixRowSep] &&3\\5&4&2&1 \end{psmatrix}} }
\def\tTrois{{\small
\begin{psmatrix}[colsep=\matrixColSep,rowsep=\matrixRowSep] &4\\5&3&2&1 \end{psmatrix}} }
\def\tDeux{{\small
\begin{psmatrix}[colsep=\matrixColSep,rowsep=\matrixRowSep] 5\\4&3&2&1 \end{psmatrix}} }

\def\tQuinze{{\small
\begin{psmatrix}[colsep=\matrixColSep,rowsep=\matrixRowSep] &&3\\5&4&1&2 \end{psmatrix}} }
\def\tTreize{{\small
\begin{psmatrix}[colsep=\matrixColSep,rowsep=\matrixRowSep] &4&2\\5&3&1 \end{psmatrix}} }
\def\tDouze{{\small
\begin{psmatrix}[colsep=\matrixColSep,rowsep=\matrixRowSep] 5&&2\\4&3&1 \end{psmatrix}} }
\def\tSept{{\small
\begin{psmatrix}[colsep=\matrixColSep,rowsep=\matrixRowSep] &4\\5&2&3&1 \end{psmatrix}} }
\def\tSix{{\small
\begin{psmatrix}[colsep=\matrixColSep,rowsep=\matrixRowSep] 5&3\\4&2&1 \end{psmatrix}} }
\def\tQuatre{{\small
\begin{psmatrix}[colsep=\matrixColSep,rowsep=\matrixRowSep] 5\\3&4&2&1 \end{psmatrix}} }

\def\tDixsept{{\small
\begin{psmatrix}[colsep=\matrixColSep,rowsep=\matrixRowSep] &4\\5&1&3&2 \end{psmatrix}} }
\def\tSeize{{\small
\begin{psmatrix}[colsep=\matrixColSep,rowsep=\matrixRowSep] 5&3\\4&1&2 \end{psmatrix}} }
\def\tQuatorze{{\small
\begin{psmatrix}[colsep=\matrixColSep,rowsep=\matrixRowSep] 5&&2\\3&4&1 \end{psmatrix}} }
\def\tVingtetun{{\small
\begin{psmatrix}[colsep=\matrixColSep,rowsep=\matrixRowSep] &4&3\\5&2&1 \end{psmatrix}} }
\def\tHuit{{\small
\begin{psmatrix}[colsep=\matrixColSep,rowsep=\matrixRowSep] 5\\2&4&3&1 \end{psmatrix}} }
\def\tNeuf{{\small
\begin{psmatrix}[colsep=\matrixColSep,rowsep=\matrixRowSep] 5&4\\3&2&1 \end{psmatrix}} }

\def\tVingttrois{{\small
\begin{psmatrix}[colsep=\matrixColSep,rowsep=\matrixRowSep] \ &4&3\\5&1&2 \end{psmatrix}} }
\def\tDixhuit{{\small
\begin{psmatrix}[colsep=\matrixColSep,rowsep=\matrixRowSep] 5\\1&4&3&2 \end{psmatrix}} }
\def\tDixneuf{{\small
\begin{psmatrix}[colsep=\matrixColSep,rowsep=\matrixRowSep] 5&4\\3&1&2 \end{psmatrix}} }
\def\tVingtdeux{{\small
\begin{psmatrix}[colsep=\matrixColSep,rowsep=\matrixRowSep] 5&&3\\2&4&1 \end{psmatrix}} }
\def\tDix{{\small
\begin{psmatrix}[colsep=\matrixColSep,rowsep=\matrixRowSep] 5&4\\2&3&1 \end{psmatrix}} }

\def\tVingtquatre{{\small
\begin{psmatrix}[colsep=\matrixColSep,rowsep=\matrixRowSep] 5&&3\\1&4&2 \end{psmatrix}} }
\def\tVingt{{\small
\begin{psmatrix}[colsep=\matrixColSep,rowsep=\matrixRowSep] 5&4\\1&3&2 \end{psmatrix}} }
\def\tVingtcinq{{\small
\begin{psmatrix}[colsep=\matrixColSep,rowsep=\matrixRowSep] 5&4\\2&1&3 \end{psmatrix}} }

\def\tVingtsix{{\small
\begin{psmatrix}[colsep=\matrixColSep,rowsep=\matrixRowSep] 5&4\\1&2&3 \end{psmatrix}} }

\def\ytUn{{\small
\begin{psmatrix}[colsep=\matrixColSep,rowsep=\matrixRowSep] 1&2&3&4&5 \end{psmatrix}} }

\def\ytDeux{{\small
\begin{psmatrix}[colsep=\matrixColSep,rowsep=\matrixRowSep] 2\\1&3&4&5 \end{psmatrix}} }

\def\ytTrois{{\small
\begin{psmatrix}[colsep=\matrixColSep,rowsep=\matrixRowSep] 4\\1&2&3&5 \end{psmatrix}} }

\def\ytQuatre{{\small
\begin{psmatrix}[colsep=\matrixColSep,rowsep=\matrixRowSep] 3\\1&2&4&5 \end{psmatrix}} }

\def\ytCinq{{\small
\begin{psmatrix}[colsep=\matrixColSep,rowsep=\matrixRowSep] 5\\1&2&3&4 \end{psmatrix}} }

\def\ytSix{{\small
\begin{psmatrix}[colsep=\matrixColSep,rowsep=\matrixRowSep] 2&4\\1&3&5 \end{psmatrix}} }

\def\ytSept{{\small
\begin{psmatrix}[colsep=\matrixColSep,rowsep=\matrixRowSep] 4&5\\1&2&3 \end{psmatrix}} }

\def\ytHuit{{\small
\begin{psmatrix}[colsep=\matrixColSep,rowsep=\matrixRowSep] 2&5\\1&3&4 \end{psmatrix}} }

\def\ytNeuf{{\small
\begin{psmatrix}[colsep=\matrixColSep,rowsep=\matrixRowSep] 3&4\\1&2&5 \end{psmatrix}} }

\def\ytDix{{\small
\begin{psmatrix}[colsep=\matrixColSep,rowsep=\matrixRowSep] 3&5\\1&2&4 \end{psmatrix}} }

\def\ytOnze{{\small
\begin{psmatrix}[colsep=\matrixColSep,rowsep=\matrixRowSep] 4\\2\\1&3&5 \end{psmatrix}} }

\def\ytDouze{{\small
\begin{psmatrix}[colsep=\matrixColSep,rowsep=\matrixRowSep] 3\\2\\1&4&5 \end{psmatrix}} }

\def\ytTreize{{\small
\begin{psmatrix}[colsep=\matrixColSep,rowsep=\matrixRowSep] 5\\2\\1&3&4 \end{psmatrix}} }

\def\ytQuatorze{{\small
\begin{psmatrix}[colsep=\matrixColSep,rowsep=\matrixRowSep] 4\\3\\1&2&5 \end{psmatrix}} }

\def\ytQuinze{{\small
\begin{psmatrix}[colsep=\matrixColSep,rowsep=\matrixRowSep] 5\\4\\1&2&3 \end{psmatrix}} }

\def\ytSeize{{\small
\begin{psmatrix}[colsep=\matrixColSep,rowsep=\matrixRowSep] 5\\3\\1&2&4 \end{psmatrix}} }

\def\ytDixSept{{\small
\begin{psmatrix}[colsep=\matrixColSep,rowsep=\matrixRowSep] 4\\2&5\\1&3 \end{psmatrix}} }

\def\ytDixHuit{{\small
\begin{psmatrix}[colsep=\matrixColSep,rowsep=\matrixRowSep] 5\\2&4\\1&3 \end{psmatrix}} }

\def\ytDixNeuf{{\small
\begin{psmatrix}[colsep=\matrixColSep,rowsep=\matrixRowSep] 4\\3&5\\1&2 \end{psmatrix}} }

\def\ytVingt{{\small
\begin{psmatrix}[colsep=\matrixColSep,rowsep=\matrixRowSep] 3\\2&5\\1&4 \end{psmatrix}} }

\def\ytVingtEtUn{{\small
\begin{psmatrix}[colsep=\matrixColSep,rowsep=\matrixRowSep] 5\\3&4\\1&2 \end{psmatrix}} }

\def\ytVingtDeux{{\small
\begin{psmatrix}[colsep=\matrixColSep,rowsep=\matrixRowSep] 4\\3\\2\\1&5 \end{psmatrix}} }

\def\ytVingtTrois{{\small
\begin{psmatrix}[colsep=\matrixColSep,rowsep=\matrixRowSep] 5\\4\\2\\1&3 \end{psmatrix}} }

\def\ytVingtQuatre{{\small
\begin{psmatrix}[colsep=\matrixColSep,rowsep=\matrixRowSep] 5\\3\\2\\1&4 \end{psmatrix}} }

\def\ytVingtCinq{{\small
\begin{psmatrix}[colsep=\matrixColSep,rowsep=\matrixRowSep] 5\\4\\3\\1&2 \end{psmatrix}} }

\def\ytVingtSix{{\small
\begin{psmatrix}[colsep=\matrixColSep,rowsep=\matrixRowSep] 5\\4\\3\\2\\1 \end{psmatrix}} }

    \newcommand{\Sn}{\mathfrak{S}_n} %
    \newcommand{\mySn}[1]{\mathfrak{S}_{#1}} %
    \newcommand{\nn}{\{1, 2, ..., n\}} %

\def\Tabvrule{\vrule width-0.4pt}       
\def\Tabhrule{\hrule \hrule height-0.4pt} 
\def\Tabstrut{\vrule height2.2ex 
                     depth0.8ex  
                     width0ex    
\relax}

\def\PasCase#1{\omit%
            $\vcenter{\hbox {\vbox to 0.4pt{}}
               \hbox{\makebox[3ex]{\Tabstrut$#1$}}}%
               \Tabvrule$}
\def\PasCasePoint{\PasCase{\cdot}}
\def\DessinCarre#1{%
    \vcenter{\hbox{}\hrule
             \hbox{\vrule\makebox[3ex]{\Tabstrut$#1$}\vrule}\Tabhrule}%
             \Tabvrule}
\def\GenRuban#1{\vcenter{\halign{&$\DessinCarre{##}$\cr#1}}\egroup}

\def\sTabvrule{\vrule width-0.4pt}
\def\sTabhrule{\hrule \hrule height-0.4pt}
\def\sTabstrut{\vrule height1.6ex depth0.6ex width0ex \relax}
\def\sDessinCarre#1{%
    \vcenter{\hbox{}\hrule
             \hbox{\vrule\makebox[2.3ex]%
                  {\sTabstrut$\scriptstyle#1$}\vrule}\sTabhrule}%
             \sTabvrule}
\def\sGenRuban#1{\vcenter{\halign{&$\sDessinCarre{##}$\cr#1}}\egroup}

\def\ruban{%
  \bgroup
  \let\ =\omit
  \let\\=\cr
  \let\.=\PasCasePoint
  \offinterlineskip
  \GenRuban}

\begin{document}
\title[On the Young-Fibonacci insertion algorithm] {On the Young-Fibonacci insertion algorithm}
\author[J. Nzeutchap]{Janvier Nzeutchap}

\address{LITIS EA 4051 (Laboratoire d'Informatique, de Traitement de l'Information et des Syst\`emes)\\
Avenue de l'Universit\'e, 76800 Saint Etienne du Rouvray, France}
\email{janvier.nzeutchap@univ-mlv.fr}
\urladdr{http://monge.univ-mlv.fr$/^{\sim}$nzeutcha}


\begin{abstract}
This work is concerned with some properties of the Young-Fibonacci
insertion algorithm and its relation with Fomin's growth diagrams.
It also investigates a relation between the combinatorics of
Young-Fibonacci tableaux and the study of Okada's algebra associated
to the Young-Fibonacci lattice. The original algorithm was
introduced by Roby and we redefine it in such a way that both the
insertion and recording tableaux of any permutation are
\emph{conveniently} interpreted as chains in the Young-Fibonacci
lattice. A property of Killpatrick's evacuation is given a simpler
proof, but this evacuation is no longer needed in making Roby's and
Fomin's constructions coincide. We provide the set of
Young-Fibonacci tableaux of size $n$ with a structure of graded
poset, induced by the weak order on permutations of the symmetric
group, and realized by transitive closure of elementary
transformations on tableaux. We show that this poset gives a
combinatorial interpretation of the coefficients in the transition
matrix from the analogue of complete symmetric functions to analogue
of the Schur functions in Okada's algebra. We end with a quite
similar observation for four posets on Young-tableaux studied by
Taskin.
\end{abstract}

\maketitle

\vspace{-1cm}

\tableofcontents

\vspace{-1.5cm}

\section{Introduction}

The Young lattice ($\youngLattice$) is defined on the set of
partitions of positive integers, with covering relations given by
the natural inclusion order. The differential poset nature of this
graph was generalized by Fomin who introduced graph duality
\cite{graph_duality}. With this extension he introduced
\cite{schensted_for_dual_graphs} a generalization of the classical
Robinson-Schensted-Knuth \cite{schensted_subsequence, knuth_rsk}
algorithm, giving a general scheme for establishing bijective
correspondences between couples of saturated chains in dual graded
graphs, both starting at a vertex of rank 0 and having a common end
point of rank $n$, on the one hand, and permutations of the
symmetric group $\Sn$ on the other hand. This approach naturally
leads to the Robinson-Schensted insertion algorithm.

Roby \cite{roby_thesis} gave an insertion algorithm analogous to the
Schensted correspondence, which maps a permutation $\sigma$ onto a
couple made of a Young-Fibonacci tableau $P(\sigma)$ and a path
tableau $Q(\sigma)$. Roby's path tableau $Q(\sigma)$ is canonically
interpreted as a saturated chain in the Fibonacci lattice $Z(1)$
introduced by Stanley \cite{fibolattice_lattice} and also by Fomin
\cite{generalized_rsk}. Roby also showed that Fomin's approach is
partially equivalent to his construction.

Indeed in Roby's construction, only the saturated chain
$\fominQFibo$ obtained from Fomin's growth diagram has an
interpretation as a representation of the path tableau $Q(\sigma)$,
while there seems to be no way to translate the Young-Fibonacci
tableau $P(\sigma)$ into its equivalent chain $\fominPFibo$.
Contrarily to the approach of Killpatrick
\cite{evacuation_fibonacci} who has used an evacuation to relate the
two constructions of Roby and Fomin, we show that with a suitable
mechanism for converting a saturated chain in the Young-Fibonacci
lattice into a Young-Fibonacci tableau, Roby's construction
naturally coincides with Fomin's one.

The paper is organized as follows. In Section
\ref{section::young::fibo::lattice} we recall the definition of the
Young-Fibonacci lattice, then in Section
\ref{section::young::fibo::tableaux} we define a mechanism for
converting a saturated chain in this lattice into a standard
Young-Fibonacci tableau. In the same section, we also introduce a
modification in Roby's algorithm, in such a way that both the
insertion and recording tableaux of any permutation will have an
interpretation in terms of saturated chain in the Young-Fibonacci
lattice. 
In Section \ref{section::growth::poset} we relate Roby's algorithm
with Fomin's construction using growth diagrams and we compare it to
Killpatrick's work. In Section \ref{section::young::fibo::numbers},
we define an analogue of Kostka numbers for Young-Fibonacci
tableaux, and we point out one of their relation with usual
Fibonacci numbers. In Section \ref{section::order::fibo::tableaux}
we define and we study some properties of a poset on Young-Fibonacci
tableaux. This poset turns out to be a model for the interpretation
as well as the computation of another analogue of Kostka numbers,
introduced by Okada \cite{okada_algebra_fibo} in an analogue of the
algebra of symmetric functions, associated to the Young-Fibonacci
lattice. We prove this result is Section
\ref{section::okada::algebra}, and in the last section of the paper
we prove a similar result relating usual Kostka numbers with four
posets on Young-tableaux studied by Taskin \cite{taskin}.
\subsection{The Young-Fibonacci lattice}\label{section::young::fibo::lattice}\

A \emph{Fibonacci diagram} or \emph{snakeshape} of size $n$ is a
column by column graphical representation of a composition of an
integer $n$, with parts equal to 1 or 2. The number of such
compositions is the $n^{th}$ Fibonacci number. A partial order is
defined on the set of all snakeshapes, in such a way to obtain an
analogue of the Young lattice of partitions of integers
($\youngLattice$). This lattice is called the Young-Fibonacci
lattice ($\youngFiboLattice$) and it was introduced by Stanley
\cite{fibolattice_lattice} and also by Fomin \cite{generalized_rsk}.
As we will see in the sequel, there is a considerable similarity
between the two lattices, as well as the combinatorics of tableaux
their induce. The covering relations in $\youngFiboLattice$ are
given below, for any snakeshape $u$.
\begin{enumerate}
  \item $u$ is covered by the snakeshape obtained by attaching a single box just in front ;
  \item $u$ is covered by the snakeshape obtained by adding a single box on top
  of its first single-boxed column, reading $u$ from left to right.
  \item if $u$ starts with a series of two-boxed columns,
  then it is covered by all snakeshapes obtained by inserting a
  single-boxed column just after any of those columns.
\end{enumerate}

The rank $|u|$ of a snakeshape $u$ is the sum of digits of the
corresponding Fibonacci word. Its length will be denoted $\ell(u)$.
Let $u$ and $v$ be two snakeshapes such that $v$ covers $u$ in
$\youngFiboLattice$, the cell added to $u$ to obtain $v$ is an
\emph{inner corner} of $v$, it is also called an \emph{outer corner}
of $u$.

\begin{remark}
Young-Fibonacci tableaux ($\youngFiboTableau$) will naturally appear
as numberings of snakeshapes, satisfying certain conditions
described in the sequel, the same way as Young tableaux are
numberings of partitions of integers with prescribed numbering
conditions. The numbering conditions of Young-Fibonacci tableaux are
deduced from the description of the Young-Fibonacci insertion
algorithm (Section \ref{section::young::fibo::algo}).
\end{remark}

Below is a pictorial representation of a finite realization of
$\youngFiboLattice$, from rank $0$ up to rank $n = 4$, with black
cells representing inner corners.

\def\colSep{0}
\def\rowSep{-0.1}
\def\corner{\blacksquare}%
\begin{figure}[h]
$$\begin{psmatrix}[colsep=0.25,rowsep=0.5]
  & [name=s1111]\snakeOneOneOneOne && [name=s211]\snakeTwoOneOne &&& [name=s22]\snakeOneTwoOne &&& [name=s31]\snakeTwoTwo &&&& [name=s4]\snakeOneOneTwo \\[0pt]
  && [name=s111]\snakeOneOneOne &&&& [name=s21]\snakeTwoOne &&&&& [name=s3]\snakeOneTwo \\[0pt]
  &&& [name=s11]\snakeOneOne &&&&&& [name=s2]\snakeTwo \\[0pt]
  &&&&&&[name=s1]\snakeOne\\[0pt]
  &&&&&& [name=s0]\emptysnake
  \psset{nodesep=5pt,arrows=-}
  \ncline{s111}{s1111} \ncline{s111}{s211} \ncline{s21}{s211}
  \ncline{s21}{s22} \ncline{s21}{s31} \ncline{s3}{s31} \ncline{s3}{s4}
  \ncline{s11}{s111} \ncline{s11}{s21} \ncline{s2}{s21} \ncline{s2}{s3}
  \ncline{s0}{s1} \ncline{s1}{s11} \ncline{s1}{s2}
\end{psmatrix}$$
\caption{The Young-Fibonacci lattice.}\end{figure}

Now let us look at the problem of converting a saturated chain in
$\youngFiboLattice$ into a standard $\youngFiboTableau$.

\section{Young-Fibonacci tableaux and Young-Fibonacci insertion
algorithm}\label{section::young::fibo::tableaux}
 In $\youngLattice$,
any saturated chain starting at the empty partition can be
canonically converted into a standard Young tableau, and this
representation is \emph{convenient} in many ways. It consists in
labeling the boxes as their occur in the chain. As already observed
by Roby \cite{roby_thesis}, one question which presents itself is to
do the same in $\youngFiboLattice$ for any saturated chain starting
at the empty snakeshape $\emptyset$. \emph{The need of such a
conversion mechanism will appear in section
\ref{section::growth::poset} in the interpretation of two saturated
chains in a growth diagram}.

One may also use the canonical labeling to convert a saturated chain
of $\youngFiboLattice$ into a tableau, but Roby had already pointed
out that one major problem with this canonical labeling is that
\emph{except for the trivial rule that each element in the top row
must be greater than the one below it, no other obvious rules govern
what numberings are allowed for a given shape}. We suggest that one
first defines simple rules governing what numberings are allowed for
a given shape, so that it be easy to decide if a numbering of a
snakeshape is a legitimate Young-Fibonacci tableau or not. The
convention we use is described in the next section.

\subsection{Converting a chain in $\youngFiboLattice$ into a
standard Young-Fibonacci tableau}\label{section::new::approach}\

Since we do not use the same conventions as Roby \cite{roby_thesis}
and Fomin \cite{graph_duality}, let us give the following definition
of Young-Fibonacci tableaux.

\begin{definition}\label{def::fibo::tableau}A numbering of a snakeshape with distinct nonnegative integers is a
standard Young-Fibonacci tableau ($\standardYoungFiboTableau$) under
the following conditions.\begin{enumerate}
  \item entries are strictly increasing in columns ;
  \item any entry on top in any column has no entry greater than itself on its right.
\end{enumerate}
\end{definition}

To convert a chain in $\youngFiboLattice$ into a standard
$\youngFiboTableau$, \emph{one will follow the canonical approach as
far as the new box added to the chain lies in the first column}.
Example with the chain $\fominQFibo = (\emptyset, 1, 2, 12, 22, 221,
2211, 21211)$ ; the sub-chain $(\emptyset, 1, 2, 12, 22)$ is
converted as follows.
{\small $$ \emptyset \,\, \to \,\, \ruban{\mbox{\bf 1}\\} \,\, \to \,\,%
\ruban{\mbox{\bf 2}\\1\\} \,\, \to \,\, \ruban{\ &2\\\mbox{\bf 3}&1\\}%
\,\, \to \,\,  \ruban{\mbox{\bf 4}&2\\3&1\\}$$}

Now moving from the shape $22$ to the shape $221$ in
$\youngFiboLattice$, one inserts a box just after a two-boxed column
of the previous shape. \emph{In such a situation, one will move the
entry on top in that column into the newly created box, and then
shift the other entries of the top row to the right. Finally, if $n$
is the largest entry in the partial tableau obtained, then label the
box on top in the first column with $(n+1)$}.

The conversion started above keeps on as follows, $x^k$ means that
writing or moving the label $x$ is the $k^{th}$ action performed during the current step.\\

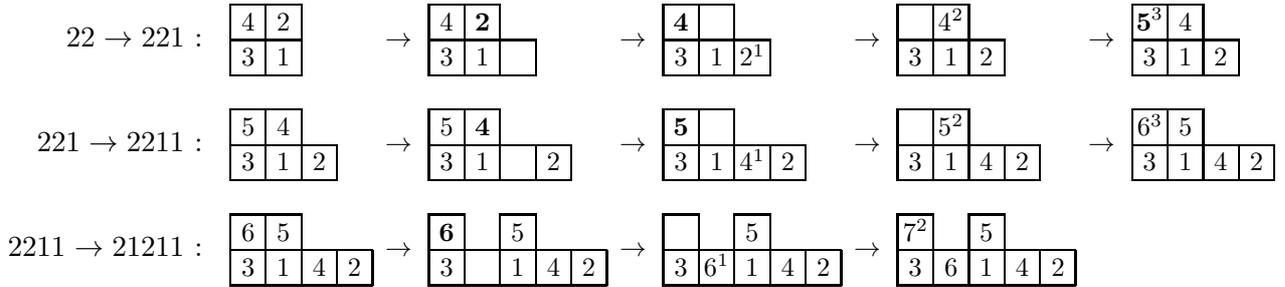
\def\negSpace{\!\!\!}
\def\mysymbol{}

\noindent\begin{tabular}{rlllll} $22 \to 221$ :&{\small
$\ruban{4&2\\3&1\\}$}&{\small $\negSpace \to \,\, \ruban{4&\mbox{\bf
2}\\3&1&\mysymbol\\}$}& {\small $\negSpace \to \,\, \ruban{\mbox{\bf
4}&\mysymbol\\3&1&2^{1}\\}$ } & {\small $\negSpace \to \,\,
\ruban{\mysymbol&4^2\\3&1&2\\}$ } &
 {\small $\negSpace \to \,\, \ruban{{\bf 5}^3&4\\3&1&2\\}$ }\\ \\
$221 \to 2211$ : & {\small $\ruban{5&4\\3&1&2\\}$} & {\small $
\negSpace \to \,\, \ruban{5&\mbox{\bf 4}\\3&1&\mysymbol&2\\}$} &
{\small $\negSpace \to \,\, \ruban{\mbox{\bf
5}&\mysymbol\\3&1&4^1&2\\}$} & {\small $\negSpace \to \,\,
\ruban{\mysymbol&5^2\\3&1&4&2\\}$} & {\small $\negSpace \to \,\,
\ruban{6^3&5\\3&1&4&2\\}$}\\ \\
$2211 \to 21211$ : &  {\small $\ruban{6&5\\3&1&4&2\\}$} &  {\small $
\negSpace \to \,\, \ruban{\mbox{\bf 6}&\ &5\\3&\mysymbol&1&4&2\\}$}
&
{\small $\negSpace \to \,\, \ruban{\mysymbol&\ &5\\3&6^1&1&4&2\\}$} & {\small $\negSpace \to \,\, \ruban{7^2&\ &5\\3&6&1&4&2\\}$}\\
\end{tabular}\\

It easily follows from the description above that this mechanism
produces only legitimate $\youngFiboTableau$ (Definition
\ref{def::fibo::tableau}), and that the conversion is reversible.
Now another question which presents itself is how to count standard
$\youngFiboTableau$ of a given shape $u \neq \emptyset$, we denote
this number by $\fiboTableauNumber{u}$. Let us first recall the
formula counting linear extensions of a binary tree poset $\posetp$.
\begin{equation}\label{equation::hook:trees}
|\linearext{\posetp}| \,\, = \,\, \frac{n\, !}{d_1 d_2 \cdots d_n}
\end{equation}
where for the $i^{th}$ node $v_i$, $d_i$ is the number of nodes $v
\leq_{\posetp} v_i$. This formula is due to Knuth
\cite{art_of_comp_programming}, and since any snakeshape $u$ can be
canonically assimilated to a poset $\posetp_u$, then we have the
following.

\begin{proposition}Standard Young-Fibonacci tableaux of a given shape
are counted by the hook-length formula for binary trees.
\end{proposition}

To apply the formula to a snakeshape $u$, count it cells from right
to left and from bottom to top, labeling the first box and each box
appearing in the bottom row of any two-boxed column. The number of
standard $\youngFiboTableau$ of the given shape is the product of
all the labels obtained.
\begin{example}\label{example::hook}Let us consider $u = 2212$.\end{example}
\vspace{-0.75cm}
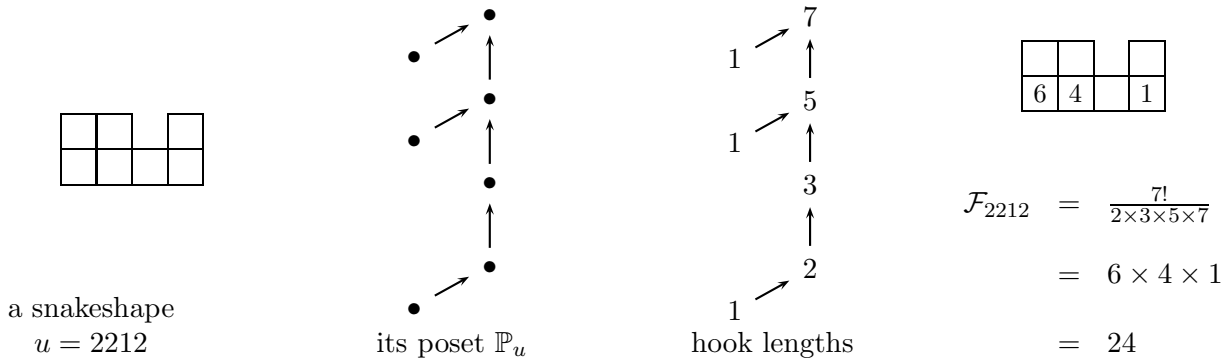
\begin{multicols}{4}
{\small
$$\ruban{&&\ &\\&&&\\}$$}
\begin{tabular}{c}
a snakeshape\\
$u = 2212$\\
\end{tabular}

\columnbreak

\def\eUn{\bullet} \def\eDeux{\bullet} \def\eTrois{\bullet} \def\eQuatre{\bullet}
\def\eCinq{\bullet} \def\eSix{\bullet} \def\eSept{\bullet}

\def\myposet{
\begin{psmatrix}[colsep=0.4,rowsep=0.1]
&& [name=e4]\eQuatre\\[0pt]
[name=e7]\eSept\\[0pt]
&& [name=e1]\eUn\\[0pt]
[name=e6]\eSix\\[0pt]
&& [name=e5]\eCinq\\[0pt]
\\[0pt]
&& [name=e2]\eDeux\\[0pt]
[name=e3]\eTrois
  \psset{nodesep=5pt,arrows=->}
  \ncline{e3}{e2} \ncline{e2}{e5} \ncline{e5}{e1} \ncline{e6}{e1} \ncline{e1}{e4} \ncline{e7}{e4}
\end{psmatrix}
}

$$\myposet$$
$$\mbox{its poset } \posetp_u$$

\columnbreak
\def\eUn{5} \def\eDeux{2} \def\eTrois{1} \def\eQuatre{7}
\def\eCinq{3} \def\eSix{1} \def\eSept{1}
$$\myposet$$
$$\mbox{hook lengths}$$

\columnbreak

$${\small \ruban{&&\ &\\6&4&&1\\}} $$

$$
\begin{array}{lcl}
\fiboTableauNumber{2212} & = & \frac{7!}{2 \times 3 \times 5 \times 7}\\
&\\
& = & 6 \times 4 \times 1\\
& \\
& = & 24
\end{array}
$$

\end{multicols}

\subsection{Redefining the Young-Fibonacci Insertion
Algorithm}\label{section::young::fibo::algo}\

We refer the reader to \cite{roby_thesis, evacuation_fibonacci} for
a description of the original algorithm ; below is the one we
consider.

\begin{definition}\label{def::young::fibo::algo}
The Young-Fibonacci insertion algorithm maps a permutation $\sigma$
onto a \emph{couple of standard} $\youngFiboTableau$ built as
follows. The insertion tableau $P(\sigma)$ is built by reading
$\sigma$ \underline{from right to left}, matching any of the letters
encountered (and not yet matched) with the maximal one (not yet
matched) on its left if any, provided that the latter be greater
than the first. The recording tableau $Q(\sigma)$ records the
positions of the letters, in the reverse order of the one in which
they are matched.
\end{definition}
\begin{example}For $\sigma = 2715643$, we have the following.\end{example}\vspace{-0.5cm}

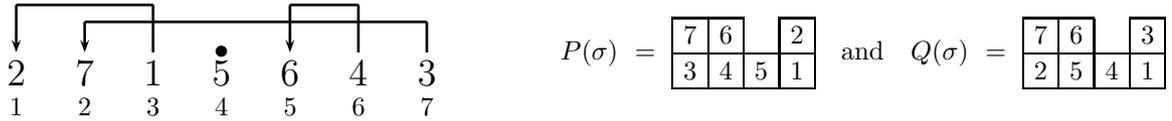
\begin{multicols}{2}
$$
\psset{unit=0.45cm}
\def\lineOne{0} \def\lineTwo{1} \def\lineThree{1.6} \def\lineFour{2.5}
\def\lineFive{3} \def\lineSix{3.5}
\begin{pspicture}(0,0)(12,\lineFive)
\def\fontUp{\Large} \def\fontDown{\small} \def\alone{$\bullet$}
\rput(0,\lineTwo){{\fontUp 2}} \rput(2,\lineTwo){{\fontUp 7}}
\rput(4,\lineTwo){{\fontUp 1}} \rput(6,\lineTwo){{\fontUp 5}}
\rput(8,\lineTwo){{\fontUp 6}} \rput(10,\lineTwo){{\fontUp 4}} \rput(12,\lineTwo){{\fontUp 3}}%
\rput(0,\lineOne){{\fontDown 1}} \rput(2,\lineOne){{\fontDown 2}}
\rput(4,\lineOne){{\fontDown 3}} \rput(6,\lineOne){{\fontDown 4}}
\rput(8,\lineOne){{\fontDown 5}} \rput(10,\lineOne){{\fontDown 6}} \rput(12,\lineOne){{\fontDown 7}}%
\psline[linewidth=1pt]{-}(12,\lineThree)(12,\lineFour)
\psline[linewidth=1pt]{-}(12,\lineFour)(2,\lineFour)
\psline[linewidth=1pt]{->}(2,\lineFour)(2,\lineThree)
\psline[linewidth=1pt]{-}(10,\lineThree)(10,\lineFive)
\psline[linewidth=1pt]{-}(10,\lineFive)(8,\lineFive)
\psline[linewidth=1pt]{->}(8,\lineFive)(8,\lineThree)
\rput(6,\lineThree){\alone}%
\psline[linewidth=1pt]{-}(4,\lineThree)(4,\lineFive)
\psline[linewidth=1pt]{-}(4,\lineFive)(0,\lineFive)
\psline[linewidth=1pt]{->}(0,\lineFive)(0,\lineThree)
\end{pspicture}
$$

\columnbreak

{\small $$ P(\sigma) \,\, = \,\, \ruban{7&6&\ &2\\3&4&5&1\\} %
\quad \mbox{and} \quad Q(\sigma) \,\, = \,\, \ruban{7&6&\
&3\\2&5&4&1\\}
$$ }
\end{multicols}

\begin{remark}That both $P(\sigma)$ and $Q(\sigma)$ are standard Young-Fibonacci tableaux
(Definition \ref{def::fibo::tableau}) is clear from the description
of the algorithm. This is not the case in the original algorithm
where $P(\sigma)$ and $Q(\sigma)$ are not of the same type. Indeed,
with the original insertion algorithm, the insertion tableau is the
same as the tableau $P(\sigma)$ above, but the recording tableau
${\bf Q}(\sigma)$ which follows does not satisfy Definition
\ref{def::fibo::tableau}. {\small
$${\bf Q_{Roby}}(\sigma) = \ruban{3&7&\ &4\\2&6&5&1\\}
$$}\end{remark}

The definition of $Q(\sigma)$ we adopt is inspired from the
hypoplactic \cite{ncsf4} and sylvester \cite{f_hivert_pbt} insertion
algorithms, where $Q(\sigma)$ also records the positions in $\sigma$
of the labels of $P(\sigma)$. With this definition, some essential
properties of the Young-Fibonacci correspondence have a much easier
combinatorial proof, which is not always the case in
\cite{roby_thesis}. For example, let us recall the involution
property.

\begin{theorem}\cite{roby_thesis}\label{proposition::schutzFibo}
For any permutation $\sigma$, $P(\sigma^{-1}) =
Q(\sigma)$.\end{theorem}

\begin{proof}Consider the geometric construction by
Killpatrick \cite{evacuation_fibonacci}, and recall that $P(\sigma)$
corresponds to reading vertical coordinates of the rightmost and
uppermost $\permOneSymb$ in that order, for any broken line. As for
$Q(\sigma)$, we have defined it in such a way that it corresponds to
reading horizontal coordinates of the uppermost and rightmost
$\permOneSymb$ in that order. The construction for $\sigma^{-1}$ is
obtained by transposing the one for $\sigma$.
\end{proof}

Another fundamental property of Roby's algorithm which is easily
proved using Definition \ref{def::young::fibo::algo} follows.


\begin{theorem}\cite{roby_thesis}\label{theorem::invol}
Let $\sigma$ be an involution of the symmetric group, then the cycle
decomposition of $\sigma$ is the column reading of its insertion
tableau $P(\sigma)$.\end{theorem}

We give two other canonical words associated with a tableau $t$ ; so
if we let $\fiboClass{t}$ denotes the equivalence class made of
permutations having $t$ as insertion tableau, then $\fiboClass{t}$
has at least three canonical elements. The first canonical element
is its canonical involution, that is the only involution the cycles
of which coincide with the columns of $t$, as stated in Theorem
\ref{theorem::invol}. The two other canonical elements are the
maximal (resp. minimal) element for the lexicographical order. We
will make use of these elements in Section
\ref{section::order::fibo::tableaux}.

\begin{lemma}\label{lemma::cano::max}Let $t$ be a Young-Fibonacci
tableau, $w_1$ the left-to-right reading of its top row and $w_2$
the \underline{right-to-left} reading of its bottom row, then
$w_1.w_2$ (where $.$ denotes the usual concatenation of words) is
the maximal element (for the lexicographical order) of
$\fiboClass{t}$, it is denoted $\maxCano{t}$.\end{lemma}

\begin{lemma}\label{lemma::cano::min}The word consisting of the right-to-left and up-down column
reading of $t$ is the minimal element (for the lexicographical
order) of $\fiboClass{t}$, it is denoted $\minCano{t}$.\end{lemma}

\begin{proof}Clear from the description of the Young-Fibonacci
insertion algorithm.\end{proof} \vskip 5pt

\noindent An example is given with the tableau $t$ below ; its
canonical involution is $(13)(26)(48)(5)(7) = 36185274$, the maximal
canonical element is $86315274$, and the minimal one is $31562784$.

{\small
$$t \,\, = \,\, \ruban{8&\ &6&\ &3\\4&7&2&5&1\\}$$}

We will see (Theorem \ref{theo::class::linear::extension}) that
$\fiboClass{t}$ is the set of linear extensions of a given poset,
and additionally, $\fiboClass{t}$ is an interval of the weak order
on the symmetric group (Theorem \ref{theorem::classes::intervals}).
\section{Young-Fibonacci insertion and growth in differential
posets}

In this section we show that with the modification we have
introduced in Roby's original insertion algorithm, together with the
conversion mechanism discussed in section
\ref{section::new::approach}, the Young-Fibonacci insertion
algorithm naturally coincides with Fomin's approach using growth
diagrams. So we claim that Killpatrick's evacuation
\cite{evacuation_fibonacci} is no longer needed in making the two
constructions coincide. We give a simplification of Killpatrick's
theorem relating Roby's original algorithm to Fomin's one through an
evacuation process, and we will later need this evacuation in the
proof of Theorem \ref{theo::poset::tableaux::kostka} giving a
combinatorial interpretation of Okada's analogue of Kostka numbers.


Let us recall that Fomin's construction with growth diagrams
consists in using some \emph{local rules} in filling a diagram
giving rise to a pair of saturated chains in $\youngLattice$. For
any permutation $\sigma$, the growth diagram $d(\sigma)$ is build
the following way. First draw the permutation matrix of $\sigma$ ;
next fill the left and lower boundary of $d(\sigma)$ with the empty
snakeshape $\emptyset$. The rest of the construction is iterative ;
$d(\sigma)$ is filled from the lower left corner to the upper right
corner, following the diagonal. At each step and for any
configuration as pictured below, $z$ is obtained by applying the
local rules to the vertices $t$, $x$, $y$ and the permutation matrix
element $\alpha$. We refer the reader to
\cite{schensted_for_dual_graphs} for more details on this
construction. {\small
$$\begin{psmatrix}[colsep=0.6,rowsep=0.2]
 && b_2 && \\[0pt]
& [name=x]x && [name=z]\fbox{z} \\[0pt]
a_1 && \alpha && b_1 \\[0pt]
& [name=t]t && [name=y]y \\[0pt]
 && a_2 &&
  \psset{nodesep=5pt,arrows=->}
  \ncline{t}{x} \ncline{x}{z}
  \ncline{t}{y} \ncline{y}{z}
\end{psmatrix}$$} \vspace{-0.5cm}

\begin{figure}[h]
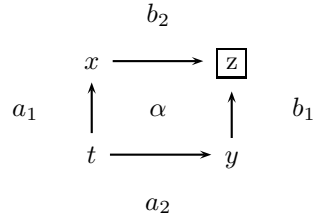

\caption{A square in a growth diagram.}\end{figure}

\begin{algorithm}
\caption{: local rules for $\youngFiboLattice$}
\label{algo::rcor::abr}
\begin{algorithmic}[1]
\IF{$x \neq y$ and $y \neq t$}%
    \STATE $z \Affect t$, with a two-boxed column added in front%
\ELSE %
    \IF{$x = y = t$ and $\alpha = 1$}%
        \STATE $z \Affect t$, with a single-boxed column added in front%
    \ELSE %
        \STATE $z$ is defined in such a way that the edge $b_i$ is degenerated whenever $a_i$ is degenerated%
    \ENDIF%
\ENDIF
\end{algorithmic}
\end{algorithm}
\subsection{Equivalence between Roby's and Fomin's
constructions}\label{section::growth::poset}\

Let us build Fomin's growth diagram for the permutation $\sigma =
2715643$.

\def\colSep{-0.04}
\def\rowSep{-0.2}
\def\corner{\square}%

{\small \centerline{
\newdimen\vcadre\vcadre=0.15cm 
\newdimen\hcadre\hcadre=0.15cm 
\setlength\unitlength{1.2mm} $\xymatrix@R=-0.5mm@C=0mm{
*{\GrTeXBox{\emptysnake}}\arx1[dd]\arx1[rr] && *{\GrTeXBox{\snakeOne}}\arx1[dd]\arx1[rr] && *{\GrTeXBox{\snakeOneOne}}\arx1[dd]\arx1[rr] &&%
*{\GrTeXBox{\snakeTwoOne}}\arx1[dd]\arx1[rr] && *{\GrTeXBox{\snakeTwoTwo}}\arx1[dd]\arx1[rr] && *{\GrTeXBox{\snakeTwoOneTwo}}\arx1[dd]\arx1[rr] &&%
*{\GrTeXBox{\snakeTwoOneOneTwo}}\arx1[dd]\arx1[rr] && *{\GrTeXBox{\snakeTwoTwoOneTwo}}\arx1[dd]\\
& *{\GrTeXBox{}} && *{\GrTeXBox{\permOneSymb}}\arx3[rrrrrrrrrr] && *{\GrTeXBox{}} &&%
*{\GrTeXBox{}} && *{\GrTeXBox{}} && *{\GrTeXBox{}} &&%
*{\GrTeXBox{}}\\
*{\GrTeXBox{\emptysnake}}\arx1[dd]\arx1[rr] && *{\GrTeXBox{\snakeOne}}\arx1[dd]\arx1[rr] && *{\GrTeXBox{\snakeOne}}\arx1[dd]\arx1[rr] &&%
*{\GrTeXBox{\snakeTwo}}\arx1[dd]\arx1[rr] && *{\GrTeXBox{\snakeOneTwo}}\arx1[dd]\arx1[rr] && *{\GrTeXBox{\snakeOneOneTwo}}\arx1[dd]\arx1[rr] &&%
*{\GrTeXBox{\snakeTwoOneTwo}}\arx1[dd]\arx1[rr] && *{\GrTeXBox{\snakeTwoTwoTwo}}\arx1[dd]\\
& *{\GrTeXBox{}} && *{\GrTeXBox{}} && *{\GrTeXBox{}} &&%
*{\GrTeXBox{}} && *{\GrTeXBox{\permOneSymb}}\arx3[rr] && *{\GrTeXBox{}} &&%
*{\GrTeXBox{}}\\
*{\GrTeXBox{\emptysnake}}\arx1[dd]\arx1[rr] && *{\GrTeXBox{\snakeOne}}\arx1[dd]\arx1[rr] && *{\GrTeXBox{\snakeOne}}\arx1[dd]\arx1[rr] &&%
*{\GrTeXBox{\snakeTwo}}\arx1[dd]\arx1[rr] && *{\GrTeXBox{\snakeOneTwo}}\arx1[dd]\arx1[rr] && *{\GrTeXBox{\snakeOneTwo}}\arx1[dd]\arx1[rr] &&%
*{\GrTeXBox{\snakeTwoTwo}}\arx1[dd]\arx1[rr] && *{\GrTeXBox{\snakeTwoOneTwo}}\arx1[dd]\\
& *{\GrTeXBox{}} && *{\GrTeXBox{}} && *{\GrTeXBox{}} &&%
*{\GrTeXBox{\permOneSymb}} && *{\GrTeXBox{}} && *{\GrTeXBox{}} &&%
*{\GrTeXBox{}}\\
*{\GrTeXBox{\emptysnake}}\arx1[dd]\arx1[rr] && *{\GrTeXBox{\snakeOne}}\arx1[dd]\arx1[rr] && *{\GrTeXBox{\snakeOne}}\arx1[dd]\arx1[rr] &&%
*{\GrTeXBox{\snakeTwo}}\arx1[dd]\arx1[rr] && *{\GrTeXBox{\snakeTwo}}\arx1[dd]\arx1[rr] && *{\GrTeXBox{\snakeTwo}}\arx1[dd]\arx1[rr] &&%
*{\GrTeXBox{\snakeOneTwo}}\arx1[dd]\arx1[rr] && *{\GrTeXBox{\snakeTwoTwo}}\arx1[dd]\\
& *{\GrTeXBox{}} && *{\GrTeXBox{}} && *{\GrTeXBox{}} &&%
*{\GrTeXBox{}} && *{\GrTeXBox{}} && *{\GrTeXBox{\permOneSymb}}\arx3[uuuu] &&%
*{\GrTeXBox{}}\\
*{\GrTeXBox{\emptysnake}}\arx1[dd]\arx1[rr] && *{\GrTeXBox{\snakeOne}}\arx1[dd]\arx1[rr] && *{\GrTeXBox{\snakeOne}}\arx1[dd]\arx1[rr] &&%
*{\GrTeXBox{\snakeTwo}}\arx1[dd]\arx1[rr] && *{\GrTeXBox{\snakeTwo}}\arx1[dd]\arx1[rr] && *{\GrTeXBox{\snakeTwo}}\arx1[dd]\arx1[rr] &&%
*{\GrTeXBox{\snakeTwo}}\arx1[dd]\arx1[rr] && *{\GrTeXBox{\snakeOneTwo}}\arx1[dd]\\
& *{\GrTeXBox{}} && *{\GrTeXBox{}} && *{\GrTeXBox{}} &&%
*{\GrTeXBox{}} && *{\GrTeXBox{}} && *{\GrTeXBox{}} &&%
*{\GrTeXBox{\permOneSymb}}\arx3[uuuuuuuu]\\
*{\GrTeXBox{\emptysnake}}\arx1[dd]\arx1[rr] && *{\GrTeXBox{\snakeOne}}\arx1[dd]\arx1[rr] && *{\GrTeXBox{\snakeOne}}\arx1[dd]\arx1[rr] &&%
*{\GrTeXBox{\snakeTwo}}\arx1[dd]\arx1[rr] && *{\GrTeXBox{\snakeTwo}}\arx1[dd]\arx1[rr] && *{\GrTeXBox{\snakeTwo}}\arx1[dd]\arx1[rr] &&%
*{\GrTeXBox{\snakeTwo}}\arx1[dd]\arx1[rr] && *{\GrTeXBox{\snakeTwo}}\arx1[dd]\\
& *{\GrTeXBox{\permOneSymb}}\arx3[rrrr] && *{\GrTeXBox{}} && *{\GrTeXBox{}} &&%
*{\GrTeXBox{}} && *{\GrTeXBox{}} && *{\GrTeXBox{}} &&%
*{\GrTeXBox{}}\\
*{\GrTeXBox{\emptysnake}}\arx1[dd]\arx1[rr] && *{\GrTeXBox{\emptysnake}}\arx1[dd]\arx1[rr] && *{\GrTeXBox{\emptysnake}}\arx1[dd]\arx1[rr] &&%
*{\GrTeXBox{\snakeOne}}\arx1[dd]\arx1[rr] && *{\GrTeXBox{\snakeOne}}\arx1[dd]\arx1[rr] && *{\GrTeXBox{\snakeOne}}\arx1[dd]\arx1[rr] &&%
*{\GrTeXBox{\snakeOne}}\arx1[dd]\arx1[rr] && *{\GrTeXBox{\snakeOne}}\arx1[dd]\\
& *{\GrTeXBox{}} && *{\GrTeXBox{}} && *{\GrTeXBox{\permOneSymb}}\arx3[uu] &&%
*{\GrTeXBox{}} && *{\GrTeXBox{}} && *{\GrTeXBox{}} &&%
*{\GrTeXBox{}}\\
*{\GrTeXBox{\emptysnake}}\arx1[rr] && *{\GrTeXBox{\emptysnake}}\arx1[rr] && *{\GrTeXBox{\emptysnake}}\arx1[rr] &&%
*{\GrTeXBox{\emptysnake}}\arx1[rr] && *{\GrTeXBox{\emptysnake}}\arx1[rr] && *{\GrTeXBox{\emptysnake}}\arx1[rr] &&%
*{\GrTeXBox{\emptysnake}}\arx1[rr] && *{\GrTeXBox{\emptysnake}}\\
 }$ } } \vspace{-0.5cm}

\begin{figure}[h]
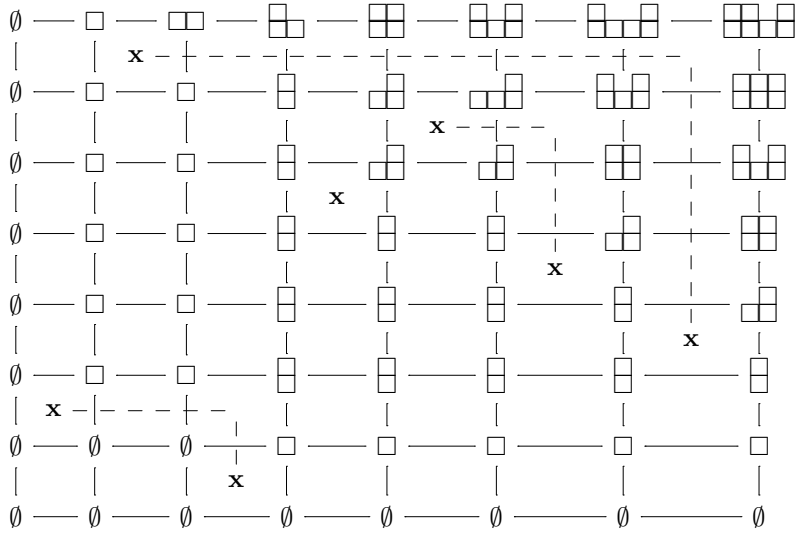

\caption{Example of growth diagram for the Young-Fibonacci
insertion.}\end{figure}

We get the paths $\fominQFibo = (\emptyset, 1, 2, 12, 22, 221, 2211,
21211)$ and $\fominPFibo = (\emptyset, 1, 11, 21, 211, 1211, 2211,
21211)$ on the upper and right boundary respectively. Now let us
convert them into Young-Fibonacci tableaux, using the mechanism
discussed in section \ref{section::new::approach}.
\def\myspace{\,}
{\small $$ \emptyset \myspace \to \myspace \ruban{1\\} \myspace \to
\myspace \ruban{2&1\\} \myspace \to \myspace \ruban{3\\2&1\\}
\myspace \to \myspace \ruban{4&3\\2&1\\} \myspace \to \myspace
\ruban{5&\ &3\\2&4&1\\} \myspace \to \myspace \ruban{6&\ &\
&3\\2&5&4&1\\} \myspace \to \myspace \ruban{7&6&\ &3\\2&5&4&1\\}
\,\, = \,\, \fominQFibo(\sigma) \,\, = \,\, Q(\sigma)
$$ }

{\small $$ \emptyset \myspace \to \myspace \ruban{1\\} \myspace \to
\myspace \ruban{2\\1\\} \myspace \to \myspace \ruban{\ &2\\3&1\\}
\myspace \to \myspace \ruban{4&2\\3&1\\} \myspace \to \myspace
\ruban{5&\ &2\\3&4&1\\} \myspace \to \myspace \ruban{6&5&2\\3&4&1\\}
\myspace \to \myspace \ruban{7&6&\ &2\\3&4&5&1\\} \,\, = \,\,
\fominPFibo(\sigma) \,\, = \,\, P(\sigma)
$$ }

So as we can see on this example, the two constructions naturally
coincide.

\begin{remark}
Let us mention that because Roby used the canonical labeling to
convert a chain into a tableau, there seemed to be no way to convert
the chain $\fominPFibo$ into its equivalent tableau $P(\sigma)$.
Killpatrick's algorithm was then an approach to relate $\fominPFibo$
with $P(\sigma)$. Our own approach consists in the introduction of a
modification of the original algorithm, and a new labeling process.
\end{remark}

\begin{theorem}Let $(\fominPFibo(\sigma), \fominQFibo(\sigma))$ be
the pair of Young-Fibonacci tableaux obtained from the permutation
$\sigma$ by using Fomin's growth diagram and let $(P(\sigma),
Q(\sigma))$ be the Young-Fibonacci insertion and recording tableaux
using Roby's insertion modified (Definition
\ref{def::young::fibo::algo}), then $ \fominPFibo(\sigma) \,\, =
\,\, P(\sigma)$ and $\fominQFibo(\sigma) \,\, = \,\, Q(\sigma)$.
\end{theorem}

\begin{proof}The equality $\fominPFibo(\sigma) = P(\sigma)$ follows from
that any snakeshape $\fominPFibo_k$ appearing in $\fominPFibo$ is
the shape of the tableau $P(\sigma_{/[1..k]})$ where
$\sigma_{/[1..k]}$ is the restriction of $\sigma$ to the interval
$[1..k]$. Indeed, the path $\fominPFibo$ is obtained applying to
$P(\sigma)$ the reverse process of the one described in section
\ref{section::new::approach}. In so doing, the cell added to
$\fominPFibo_k$ to get $\fominPFibo_{k+1}$ lies in the first column
when either $\sigma_{/[1..k+1]}$ ends with the letter $k+1$ or
$\sigma_{/[1..k+1]}$ does not end with the letter $k+1$ but
$\sigma_{/[1..k]}$ ends with the letter $k$. A quite similar
reasoning is used to prove the equality $\fominQFibo(\sigma) =
Q(\sigma)$.
\end{proof}
\subsection{Another viewpoint of Killpatrick's evacuation for
Young-Fibonacci tableaux}\label{section::discussion::evacuation}\

For a tableau $t$, this operation is defined only for top entries of
the columns of $t$. Let $a_0$ be such an entry, the tableau
resulting from the evacuation of $a_0$ is denoted
$\evacuatedTableauFromLetter{t}{a_0}$ and is built as follows.
\begin{enumerate}
  \item if $a_0$ is a single-boxed column, then just delete this column
  and, if this is necessary,
  shift one component of the remaining tableau to connect it with the other one (e.g of line 3 in the table below) ;
  \item otherwise, the box containing $a_0$ is emptied and one
  compares the entry $a_1$ that was just below $a_0$ with the entry
  $a_2$ on top of the column just to the right if any. If $a_2 <
  a_1$ then this terminates the evacuation process (e.g of line 4 in the table below).
  Otherwise, move $a_2$ on top of $a_1$, creating a new empty box in
  the tableau. If the new empty box is a single-boxed column, then
  this terminates the evacuation process (e.g of line 7, step 4, in the table below), otherwise, iteratively repeat the
  process with the entries just below and to the right of this new empty box.
\end{enumerate}

Let $t$ be a tableau of size $n$ and shape $u$. If one successively
evacuates the entries $n, (n-1), \cdots, 1$ from $t$, labeling the
boxes of $u$ according to the positions of the empty cells at the
end of the evacuation of entries, one gets a \emph{path tableau}
denoted $\evacuatedTableau{t}$. Recall that a path tableau is the
canonical labeling of a saturated chain.

\begin{remark}$\evacuatedTableau{t}$ is the same tableau as the one
described by Killpatrick \cite{evacuation_fibonacci}, with
Young-Fibonacci tableaux defined as in Definition
\ref{def::fibo::tableau}.
\end{remark}

\def\myemptycell{\bullet}
{\small
\begin{center}\begin{tabular}{|c|llclclcr|c|}\hline&&&&&&&&&\\
7 & $t \,\, = \,\, $ & {\small $\ruban{7&6& \ &2\\3&4&5&1\\}$} & $\to$ & %
{\small $\ruban{\myemptycell&{\bf \underline{6}}& \ &2\\{\bf 3}&4&5&1\\}$} & $\to$ & %
{\small $\ruban{6&\myemptycell& \ &2\\3&{\bf 4}&{\bf \underline{5}}&1\\}$} & $\to$ & %
{\small $\ruban{6&5& \ &2\\3&4&\myemptycell&1\\}$} & {\small $\ruban{&& \ &\\&&7&\\}$} \\&&&&&&&&&\\
6 & & {\small $\ruban{6&5&2\\3&4&1\\}$} & $\to$ & %
{\small $\ruban{\myemptycell&{\bf \underline{5}}&2\\{\bf 3}&4&1\\}$} & $\to$ & %
{\small $\ruban{5&\myemptycell&{\bf 2}\\3&{\bf 4}&1\\}$} & & & {\small $\ruban{&6& \ &\\&&7&\\}$} \\&&&&&&&&&\\
5 & & {\small $\ruban{5& \ &2\\3&4&1\\}$} & $\to$ & %
{\small $\ruban{\myemptycell& \ &2\\{\bf 3}&{\bf \underline{4}}&1\\}$} & $\to$ & %
{\small $\ruban{4& \ &2\\3&\myemptycell&1\\}$} & & & {\small $\ruban{&6& \ &\\&5&7&\\}$} \\&&&&&&&&&\\
4 & & {\small $\ruban{4&2\\3&1\\}$} & $\to$ & %
{\small $\ruban{\myemptycell&{\bf 2}\\{\bf 3}&1\\}$}&&&&& {\small $\ruban{4&6& \ &\\&5&7&\\}$} \\&&&&&&&&&\\
3 & & {\small $\ruban{\ &2\\3&1\\}$} & $\to$ & %
{\small $\ruban{\ &2\\\myemptycell&1\\}$}&&&&& {\small $\ruban{4&6& \ &\\3&5&7&\\}$} \\&&&&&&&&&\\
2 & & {\small $\ruban{2\\1\\}$} & $\to$ & %
{\small $\ruban{\myemptycell\\1\\}$} & & & && {\small $\ruban{4&6& \ &2\\3&5&7&\\}$} \\&&&&&&&&&\\
1 & & {\small $\ruban{1\\}$} & $\to$ & $\myemptycell$ &&&& $ev(t) \,\, =$ & {\small $\ruban{4&6& \ &2\\3&5&7&1\\}$} \\&&&&&&&&&\\
\hline
%
\end{tabular}\end{center}
\begin{table}[h]\label{table::young::fibo::evacuation}
\caption{{\small Evacuation on Young-Fibonacci tableaux.}}
\end{table}}

\begin{lemma}\label{lemma::avacuation}Let $w$ be a word with no letter repeated, let $a_0$ be one of its
letters appearing as a top element in a column of $P(w)$, and let
$w_0$ be the word obtained from $w$ by deleting the only occurrence
of $a_0$, then $\evacuatedTableauFromLetter{P(w)}{a_0} \,\, = \,\,
P(w_0)$.
\end{lemma}

\begin{proof}Easily from the description of the evacuation and the description of the Young-Fibonacci insertion algorithm
(Definition \ref{def::young::fibo::algo}).
\end{proof}

We give a simpler proof of the following theorem by Killpatrick,
relating $\evacuatedTableau{P(\sigma)}$ with $\fominPFibo$. Indeed,
\underline{using the canonical labeling}, Roby has converted the
path $\fominPFibo$ into a \underline{path tableau}
$\fominPFibo(\sigma)$ and

\begin{theorem}\cite{evacuation_fibonacci}
$\evacuatedTableau{P(\sigma)} \,\, = \,\, \fominPFibo(\sigma)$.
\end{theorem}

\begin{proof} Follows from Lemma \ref{lemma::avacuation} and the remark
that any snakeshape $\fominPFibo_k$ appearing in $\fominPFibo$ is
the shape of the tableau $P(\sigma_{/[1..k]})$ where
$\sigma_{/[1..k]}$ is the restriction of $\sigma$ to the interval
$[1..k]$.\end{proof}

\section{Fibonacci numbers and a statistic on Young-Fibonacci tableaux}\label{section::young::fibo::numbers}

\newcommand{\kostkanumber}[2]{K_{#1\!,\, #2}}
\newcommand{\youngfibonumber}[2]{\mathcal{N}_{#1\!,\, #2}}

In this section we point out a property of Young-Fibonacci numbers
defined as an analogue of Kostka numbers. Recall that the usual
Kostka numbers $\kostkanumber{\lambda}{\mu}$ are defined for two
partitions $\lambda$ and $\mu$ of the same integer $n$ and they
appear when expressing Schur functions $s_\lambda$ in terms of the
monomial symmetric functions $m_\mu$, and in the expression of the
complete symmetric functions $h_\mu$ in terms of Schur functions
$s_\lambda$.
\begin{equation}
s_\lambda \,\, = \,\, \sum_\mu \kostkanumber{\lambda}{\mu} \, m_\mu%
\qquad ; \qquad %
h_\mu \,\, = \,\, \sum_\mu \kostkanumber{\lambda}{\mu} \, s_\lambda
\end{equation}

We will not focus on the algebraic interpretation of the
$\kostkanumber{\lambda}{\mu}$ but rather on their combinatorial
interpretation in terms of tableaux. Indeed,
$\kostkanumber{\lambda}{\mu}$ counts the number of distinct
semi-standard Young-tableaux of shape $\lambda$ and valuation $\mu$,
that is to say with $\mu_i$ entries $i$ for $i = 1\,..\,\ell(\mu)$.
It is then natural to introduce the same definition with
Young-Fibonacci tableaux.

\begin{definition}A semi-standard Young-Fibonacci tableau is a numbering of a snakeshape with
nonnegative integers, not necessarily distinct, preserving the
conditions stated in Definition \ref{def::fibo::tableau}.
\end{definition}

\begin{definition}Let $u$ and $v$ be two snakeshapes of size $n$,
the Young-Fibonacci number associated to $u$ and $v$, denoted
$\youngfibonumber{u}{v}$ is the number of distinct semi-standard
Young-Fibonacci tableaux of shape $u$ and valuation $v$, that is to
say with $v_i$ entries $i$ for $i = 1\,..\,\ell(v)$.
\end{definition}

\noindent For example, for $u = 221$ and $v = 1211$, there are 4
distinct semi-standard Young-Fibonacci tableaux of shape $u$ and
valuation $v$. So $\youngfibonumber{221}{1211} = 4$.
$${\small \ruban{4&3\\1&2&2\\} \qquad \ruban{4&3\\2&1&2\\}
\qquad \ruban{4&3\\2&2&1\\} \qquad \ruban{4&2\\3&1&2\\}}
$$

\begin{proposition}Young-Fibonacci numbers are generated by the recurrence formulas below,
where both $u$ and $v$ are snakeshapes.
\begin{equation}\label{recurrence::my::kostka::fibo}
 \left\{
\begin{array}{lclclcl}
\youngfibonumber{\emptyset}{\emptyset} & = & 1 & ; & \youngfibonumber{2}{2} & = & 0\\
\youngfibonumber{1u}{v1} & = & \youngfibonumber{u}{v} & ; &
\youngfibonumber{1u}{v2} & = & \youngfibonumber{u}{v1}\\
\youngfibonumber{2u}{v1} & = & \sum_{w \, \in \, v^{1\!^-}} \,
\youngfibonumber{u}{w} & ; & \youngfibonumber{2u}{v2} & = & \sum_{w
\, \in \, v^{1\!^-}} \,
\youngfibonumber{u}{w1}\\
\end{array}\right.
\end{equation}
\noindent where $v^{1\!^-}$ denotes the \underline{multiset} of
snakeshapes obtained from $v$ either by deleting a single occurrence
of 1, or by decreasing a single entry not equal to 1, for example
$2112^{1\!^-} = [1112,\, 212,\, 212,\, 2111]$.
\end{proposition}

\begin{proof} Easily from the definition of Young-Fibonacci tableaux and Young-Fibonacci numbers.\end{proof}

{\small
\begin{center}\begin{tabular}{rl|ccccccccccccc}
& v = & 222 & 2211 & 2121 & 2112 & 21111 & 1221 & 1212 & 12111 &
1122 &
11211 & 11121 & 11112 & 111111\\
\hline\\
u = & 222   & 2 & 3 & 4 & 5 & 6 & 4 & 5 & 6 & 5 & 7 & 8 & 12 & {\bf 15}\\
& 2211 & 4 & 5 & 5 & 7 & 9 & 5 & 7 & 9 & 7 & 9 & 9 & 12 & {\bf 15}\\
& 2121 & 2 & 3 & 4 & 4 & 5 & 4 & 4 & 5 & 4 & 6 & 8 & 8 & {\bf 10}\\
& 2112 & 1 & 1 & 1 & 1 & 1 & 2 & 2 & 3 & 3 & 4 & 4 & 4 & {\bf 5}\\
& 21111 & 2 & 3 & 3 & 3 & 4 & 3 & 3 & 4 & 3 & 4 & 4 & 4 & {\bf 5}\\
& 1221 & 2 & 2 & 3 & 4 & 4 & 3 & 4 & 4 & 4 & 4 & 6 & 8 & {\bf 8}\\
& 1212 & 1 & 1 & 1 & 1 & 1 & 1 & 2 & 2 & 3 & 3 & 3 & 4 & {\bf 4}\\
& 12111 & 2 & 2 & 2 & 3 & 3 & 2 & 3 & 3 & 3 & 3 & 3 & 4 & {\bf 4}\\
& 1122 & 1 & 1 & 1 & 1 & 1 & 1 & 1 & 1 & 2 & 2 & 3 & 3 & {\bf 3}\\
& 11211 & 1 & 1 & 2 & 2 & 2 & 2 & 2 & 2 & 2 & 2 & 3 & 3 & {\bf 3}\\
& 11121 & 1 & 1 & 1 & 1 & 1 & 1 & 1 & 1 & 2 & 2 & 2 & 2 & {\bf 2}\\
& 11112 & \fbox{0} & \fbox{0} & \fbox{0} & \fbox{0} & \fbox{0} & 1 & 1 & 1 & 1 & 1 & 1 & 1 & {\bf 1}\\
& 111111 & {\bf 1} & {\bf 1} & {\bf 1} & {\bf 1} & {\bf 1} & {\bf 1} & {\bf 1} & {\bf 1} & {\bf 1} & {\bf 1} & {\bf 1} & {\bf 1} & {\bf 1}\\
\end{tabular}\end{center}
\begin{table}[h]\label{table::young::fibo::numbers::matrix}
\caption{{\small Matrix of Young-Fibonacci numbers for $n=6$.}}
\end{table}}

\begin{theorem}Let $n \geq 2$ be a positive integer, then the number
of couples $(u, v)$ of snakeshapes of size $n$ such that
$\youngfibonumber{u}{v} = 0$ is the $(n-2)^{th}$ Fibonacci number.
\end{theorem}

\begin{proof}The proof is done by induction on $n$.
Indeed using equation (\ref{recurrence::my::kostka::fibo}) it is
easy to see that $\youngfibonumber{u}{v} \neq 0$ whenever $u \neq
1^{n-2}2$. So the problem is equivalent to counting the number of
snakeshapes $v$ such that $\youngfibonumber{1^{^{n\!-\!2}}2}{v} =
0$. But $\youngfibonumber{1^{^{n\!-\!2}}2}{v} = 0$ if and only if
there exists a snakeshape $w$ such that $v = 2w$. Then the problem
is finally equivalent to counting the snakeshapes of size $(n-2)$,
and hence the result.
\end{proof}

\section{A weak order on Young-Fibonacci tableaux}\label{section::order::fibo::tableaux}

\def\noninversion{non-inversion}

In what follows, we introduce a partial and \emph{graded} order
denoted $\younFiboTableauOrder$ on the set $\youngFiboTab{n}$ of
Young-Fibonacci tableaux of size $n$. We will see (Theorem
\ref{theo::weakest::order}) that this partial ordering on
$\youngFiboTab{n}$ is such that the map from the \emph{weak order}
on the symmetric group $\Sn$ which sends each permutation $\sigma$
onto its Young-Fibonacci insertion tableau $P(\sigma)$ is
order-preserving. More particularly, standard Young-Fibonacci
classes on $\Sn$ are intervals of the weak order on $\Sn$. Recall
that the weak order on permutations of $\Sn$ is the transitive
closure of the relation $\sigma \permutohedronOrder \tau$ if $\tau =
\sigma \transposition{i}$ for some $i$, where $\transposition{i}$ is
the adjacent transposition $(i\,\, i\!+\!1)$. An inversion of a
permutation $\sigma$ is a couple $(j,i), 1 \leq i < j \leq n$ such
that $\sigma^{-1}(i) > \sigma^{-1}(j)$, that is to say $j$ appears
on the left of $i$ in $\sigma$. \underline{\emph{Note that this is
not the definition commonly used}}. The set of inversions of a
permutation $\sigma$ will be denoted $\invSet{\sigma}$, and the
number of inversions denoted $\invNumber{\sigma}$. We will be making
use of an analogous notion of \noninversion\ of a permutation
$\sigma$ which is a couple $(i,j), 1 \leq i < j \leq n$ such that
$\sigma^{-1}(i) < \sigma^{-1}(j)$, that is to say $i$ appears on the
left of $j$ in $\sigma$. The set of \noninversion s of a
permutation $\sigma$ will be denoted $\ordSet{\sigma}$.\\

\def\mya{{\bf a}}
\def\myb{{\bf b}}
\def\myc{{\bf c}}
\begin{definition}
To introduce $\younFiboTableauOrder$, we define the operation of
\emph{shifting an entry in a tableau $t$} as follows.
\def\myspace{\quad}
\begin{enumerate}
  \item the bottom entry \mya\ of any column of $t$ may move and bump up the entry \myc\ on its left if
  \myc\ is a single-boxed column of $t$. In the example below, the letter $1$ is the one being shifted.
  $$ {\small \ruban{5\\2&4&3&1\\}
  \myspace
  \begin{CD} @>\quad \mbox{shift the entry 1}\quad>> \end{CD}
  \myspace
  \small \ruban{5& \ &3\\2&4&1\\}}$$
  \item In the case \mya\ was the bottom entry in a two-boxed column,
  the top entry \myb\ will just fall down.
  \noindent In the two examples below, the letter $2$ (resp. $3$) is the one being shifted.
  $$ {\small \ruban{\ &4\\5&2&3&1\\} \myspace
  \begin{CD} @>\quad \mbox{shift 2}\quad>> \end{CD}
  \myspace \small \ruban{5\\2&4&3&1\\} %
  \qquad \mbox{and} \qquad  %
  \small \ruban{\ &4&2\\5&3&1\\} \myspace
  \begin{CD} @>\quad \mbox{shift 3}\quad>> \end{CD}
  \myspace \small \ruban{5&\ &2\\3&4&1\\}} %
    $$
   \item In the case the column just to the left of $\mya$ is two-boxed, with bottom entry $\myc$ and $\mya < \myc$,
   then $\mya$ may replace $\myc$ which on its turn is shifted to the right in such a
   way that if $\myc < \myb$ then $\myc$ will just replace $\mya$ ; otherwise
   $\myc$ is placed as a new single-boxed column between $\mya$ and $\myb$,
   and $\myb$ just falls down.
   In the two examples below, the letter $1$ (resp. $2$) is shifted.
   $${\small \ruban{5&4\\2&1&3\\} \myspace
   \begin{CD} @>\quad \mbox{shift 1}\quad>> \end{CD}
   \myspace \small \ruban{5&4\\1&2&3\\}
    \qquad \mbox{and} \qquad  %
   {\small \ruban{5&3\\4&2&1\\}} \myspace
   \begin{CD} @>\quad \mbox{shift 2}\quad>> \end{CD}
   \myspace \small \ruban{5\\2&4&3&1\\}} %
    $$
\end{enumerate}\
\end{definition}

\begin{remark}It easily follows from the definition that shifting an entry in a
tableau always produces a legitimate tableau of the same size. In an
analogous way, given a tableau $t$, one defines the reverse
operation of finding all the tableaux $t'$ such that shifting an
entry in $t'$ gives back $t$. For example, one will check that
$$ {\small \ruban{5& \ &2\\3&4&1\\}} \quad \mbox{is obtained from} \quad {\small \ruban{\ &4&2\\5&3&1\\}}%
\quad , \quad {\small \ruban{5& \ &2\\4&3&1\\}} \quad \mbox{and}
\quad {\small \ruban{5\\3&4&2&1\\}} \quad \mbox{by shifting 3 or
1}.$$ Finally it is clear that this operation is antisymmetric, that
is to say if $t'$ is obtained from $t$ by shifting a given entry,
then $t$ cannot be obtained from $t'$ by shifting an entry.
\end{remark}

The latter observation is enforced by the following lemma which also
defines the graduation of the poset $(\youngFiboTab{n},
\younFiboTableauOrder)$ we will soon introduce.

\begin{lemma}\label{lemma::inversions::number}Let $t_2$ be a tableau obtained by shifting an entry in a tableau
$t_1$, and let $\sigma_1$ (resp. $\sigma_2$) be the minimal
permutation canonically associated to $t_1$ (resp. $t_2$) as stated
in Lemma \ref{lemma::cano::min}, then the inversions sets of
$\sigma_1$ and $\sigma_2$ are related by the relation
$\invNumber{\sigma_2} = \invNumber{\sigma_1} + 1$.
\end{lemma}

\def\anything{\star}
\def\firstpartial{\mathcal{T}_1}
\def\secondpartial{\mathcal{T}_2}
\def\firstcano{w_1}
\def\secondcano{w_2}

\def\matrixColSep{0.1}
\def\matrixRowSep{-0.1}

\proof The proof takes into account all the situations one can
encounter in shifting an entry in $t_1$.
\begin{enumerate}
  \item $t_1 = \secondpartial\, \small{\begin{psmatrix}[colsep=\matrixColSep,rowsep=\matrixRowSep] & \anything\\c&a\end{psmatrix}} \, \firstpartial$
  and $t_2 = \secondpartial\, \small{\begin{psmatrix}[colsep=\matrixColSep,rowsep=\matrixRowSep] c & \\a&\anything\end{psmatrix}} \, \firstpartial$, where
  $\firstpartial$ and $\secondpartial$ are partial $\youngFiboTableau$
  having minimal canonical words $\firstcano$ and $\secondcano$ (see Lemma \ref{lemma::cano::min} for the definition),
  and $\anything$ means any entry preserving the conditions of Definition \ref{def::fibo::tableau}, and possibly no entry. The
  minimal permutations associated to $t_1$ and $t_2$ are $\sigma_1 = \firstcano\! \anything\! ac \secondcano$
  and $\sigma_2 = \firstcano\! \anything\! ca \secondcano$ respectively, and clearly
  $\sigma_2$ has one more inversion than $\sigma_1$.
  \vskip 5pt

  \item $t_1 = \secondpartial\, \small{\begin{psmatrix}[colsep=\matrixColSep,rowsep=\matrixRowSep] d & \anything\\c&a\end{psmatrix}} \, \firstpartial$
  and $t_2 = \secondpartial\, \small{\begin{psmatrix}[colsep=\matrixColSep,rowsep=\matrixRowSep] d\\a&c& \anything\end{psmatrix}} \, \firstpartial$,
  with $a < c < d$ ;
  one has $\sigma_1 = \firstcano\! \anything\! adc \secondcano$
  and $\sigma_2 = \firstcano\! \anything\! cda \secondcano$. The inversion $(d
  c)$ appears in $\sigma_1$ but not in $\sigma_2$, whereas the inversions $(d
  a)$ and $(c a)$ appear in $\sigma_2$ but not in $\sigma_1$ ; so $\sigma_2$ has one more inversion.
  \vskip 5pt

  \item $t_1 = \secondpartial\, \small{\begin{psmatrix}[colsep=\matrixColSep,rowsep=\matrixRowSep] d&b\\c&a\end{psmatrix}} \, \firstpartial$
  and $t_2 = \secondpartial\, \small{\begin{psmatrix}[colsep=\matrixColSep,rowsep=\matrixRowSep] d&b\\a&c\end{psmatrix}} \, \firstpartial$,
  with $a < c < b < d $ ;
  one has $\sigma_1 = \firstcano badc \secondcano$
  and $\sigma_2 = \firstcano bcda \secondcano$. The inversion $(d
  c)$ appears in $\sigma_1$ but not in $\sigma_2$, whereas the inversions $(d
  a)$ and $(c a)$ appear in $\sigma_2$ but not in $\sigma_1$ ; so $\sigma_2$ has one more inversion.
\endproof
\end{enumerate}\

\noindent We are now in position to provide $\youngFiboTab{n}$ with
a structure of poset.\vskip 10pt

\begin{definition}[weak order on $\youngFiboTab{n}$] Let $t$ and $t'$ be two tableaux of size $n$,
then $t$ is said smaller than $t'$ and we write $t
\younFiboTableauOrder t'$ if one can find a sequence $t_0 = t, t_1,
\cdots, t_k = t'$ of tableaux of size $n$ such that $t_{i+1}$ be
obtained from $t_i$ by shifting an entry, for $i$ from $0$ to $k-1$.
\end{definition}

\begin{proposition}\label{prop::graded::poset}$(\youngFiboTab{n}, \younFiboTableauOrder)$ is a
graded poset, the rank of a Young-Fibonacci tableau being the number
of inversions of its minimal canonical permutation.
\end{proposition}

\proof Follows from Lemma \ref{lemma::inversions::number}.
\endproof

\begin{remark}
Note that this remarkable property of graduation of the poset of
standard Young-Fibonacci tableaux of size $n$ does not apply to the
similar poset $\youngTab{n}$ of standard Young tableaux of size $n$.
The reader interested may refer to \cite{taskin} where Taskin
studied many nice properties of four partial orders on
$\youngTab{n}$.
\end{remark}

\def\matrixColSep{0}
\def\matrixRowSep{0}

{\small \centerline{
\newdimen\vcadre\vcadre=0.5cm 
\newdimen\hcadre\hcadre=0.5cm 
\setlength\unitlength{0.75mm}
$\xymatrix@R=1.3cm@C=0.8cm{%
 &  &  &  & *{\tVingtsix}& & & & & & \rho = 6 \\
 &  & *{\tVingtquatre}\arx1[rru]&  & *{\tVingt}\arx1[u]&  & *{\tVingtcinq}\arx1[llu]& & & & \rho = 5 \\
*{\tVingttrois}\arx1[rru]&  & *{\tDixhuit}\arx1[rru]\arx1[u]&  &
*{\tDixneuf}\arx1[u]& & *{\tVingtdeux}\arx1[u]&  &
*{\tDix}\arx1[llu]& & \rho = 4 \\
*{\tDixsept}\arx1[rru]\arx1[u]&  & *{\tSeize}\arx1[u]&  &
*{\tQuatorze}\arx1[u]& & *{\tVingtetun}\arx1[u]\arx1[llllllu]&  &
*{\tHuit}\arx1[u]\arx1[llu]& *{\tNeuf}\arx1[lu]\arx1[lllllu]& \rho = 3\\
*{\tQuinze}\arx1[rru]\arx1[u]&  & *{\tTreize}\arx1[rru]\arx1[llu]& &
*{\tDouze}\arx1[u]\arx1[llu]& & *{\tSept}\arx1[rru]\arx1[u]&  &
*{\tSix}\arx1[u]\arx1[llllllu]& *{\tQuatre}\arx1[u]\arx1[lllllu]&
\rho = 2\\
 & *{\tOnze}\arx1[rrru]\arx1[ru]\arx1[lu]&  &  & *{\tCinq}\arx1[rrrru]\arx1[rru]\arx1[llllu]&
 & *{\tTrois}\arx1[rrru]\arx1[u]\arx1[llllu]&  & *{\tDeux}\arx1[ru]\arx1[u]\arx1[llllu]& & \rho = 1 \\
 &  &  &  & *{\tUn}\arx1[rrrru]\arx1[rru]\arx1[u]\arx1[lllu]& & & & & & \rho = 0 \\
}$ }}


\begin{figure}[h]
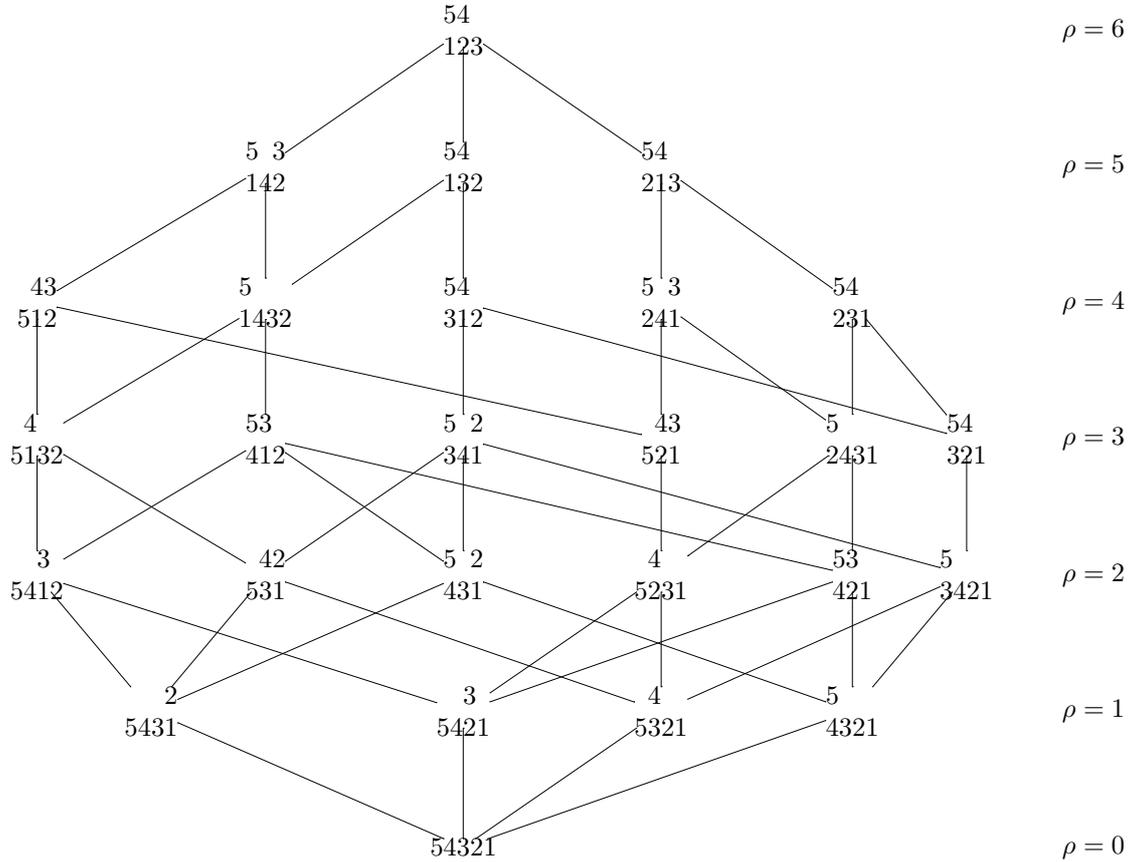

\caption{The graded weak order on Young-Fibonacci tableaux of size 5.}
\end{figure}

\begin{remark}As one will easily check it on the figure above, $(\youngFiboTab{n},
\younFiboTableauOrder)$ is not a lattice for $n = 5$ for example.
Indeed let $a = \tCinq$ and $b = \tQuatre$, then $a$ and $b$ do not
have a least upper bound.
\end{remark}

\begin{theorem}\label{theo::weakest::order}Let $t_1$ and $t_2$ be two
tableaux, then $t_1 \younFiboTableauOrder t_2$ if and only if one
can find two permutations $\tau_1$ and $\tau_2$ such that $P(\tau_1)
= t_1$, $P(\tau_2) = t_2$ and $\tau_1 \permutohedronOrder \tau_2$.
\end{theorem}

\proof It is enough to prove this statement for the case $t_2$ is
 obtained by shifting an entry in $t_1$,
 and the proof is carried out as a parallel process of the proof of Lemma
 \ref{lemma::inversions::number}. So go back to the latter proof and
 \begin{enumerate}
   \item take $\tau_i = \sigma_i$ ;
   \item take $\tau_1 = \firstcano\! \anything\! dac \secondcano$ and $\tau_2 = \firstcano\! \anything\! dca
   \secondcano$ ;
   \item take $\tau_1 = \firstcano bdac \secondcano$ and $\tau_2 = \firstcano bdca
   \secondcano$.
 \end{enumerate}
This shows that one can find two permutations $\tau_1$ and $\tau_2$
such that $P(\tau_1) = t_1$, $P(\tau_2) = t_2$ and $\tau_2 = \tau_1
\delta_i$ for some $i$, whenever $t_1 \younFiboTableauOrder
t_2$. 
Reciprocally let $\tau_1$ and $\tau_2$ be two permutations such that
$P(\tau_1) = t_1$ and $P(\tau_2) = t_2$ and $\tau_2 = \tau_1
\delta_i$ for some $i$. Then $t_2$ is obtained from $t_1$ by
shifting the entry $i$ in $t_1$.
\endproof

We now look at the structure of the Young-Fibonacci classes ; below
are two pictures of the poset $(\youngFiboTab{4},
\younFiboTableauOrder)$. On the picture on the left, vertices are
Young-Fibonacci classes corresponding to Young-Fibonacci tableaux in
the picture on the right. Recall that the rank of a class is the
number of inversions of its minimal element in the lexicographical
order. The unique involution of any class is enclosed in a
rectangle. A double edge means that there are two couples $(\tau_1,
\tau_2)$ and $(\tau{'}_1, \tau{'}_2)$ satisfying the conditions of
Theorem \ref{theo::weakest::order}.

\def\myvspace{\quad}
\def\matrixColSep{0.2}
\def\vertexRowSep{0cm}
\def\vertexColSep{0mm}

\def\rowSep{0.5cm}
\def\colSep{15mm}

\def\myposet{ {\small \centerline{
\newdimen\vcadre\vcadre=0.2cm 
\newdimen\hcadre\hcadre=0.2cm 
$\xymatrix@R=\rowSep@C=\colSep{ &
*{\GrTeXBox{\yfClassDix}}\\
*{\GrTeXBox{\yfClassHuit}} \arx2[ur] && *{\GrTeXBox{\yfClassNeuf}} \arx2[ul]\\
*{\GrTeXBox{\yfClassSept}} \arx1[u] &%
*{\GrTeXBox{\yfClassSix}} \arx2[ul] &%
*{\GrTeXBox{\yfClassQuatre}} \arx2[u]\\
*{\GrTeXBox{\yfClassCinq}} \arx1[u] \arx1[ur] &%
*{\GrTeXBox{\yfClassTrois}} \arx1[ul] \arx1[ur] &%
*{\GrTeXBox{\yfClassDeux}} \arx2[u] \arx2[ul]\\
& *{\GrTeXBox{\yfClassUn}}\arx1[ul] \arx1[u] \arx1[ur]
 }$ }} }
\begin{multicols}{2}
\myposet \columnbreak

\def\yfClassUn{{\small \ruban{4&3&2&1\\}}}
\def\yfClassDeux{{\small \ruban{4\\3&2&1\\}}}
\def\yfClassTrois{{\small \ruban{\ &3\\4&2&1\\}}}
\def\yfClassCinq{{\small \ruban{\ & \ &2\\4&3&1\\}}}
\def\yfClassQuatre{{\small \ruban{4\\2&3&1\\}}}
\def\yfClassSix{{\small \ruban{4&2\\3&1\\}}}
\def\yfClassSept{{\small \ruban{\ &3\\4&1&2\\}}}
\def\yfClassNeuf{{\small \ruban{4&3\\2&1\\}}}
\def\yfClassHuit{{\small \ruban{4\\1&3&2\\}}}
\def\yfClassDix{{\small \ruban{4&3\\1&2\\}}}

\def\rowSep{0.5cm}
\def\colSep{10mm}
\myposet

\end{multicols}

\vspace{-0.5cm}
\begin{figure}[h]
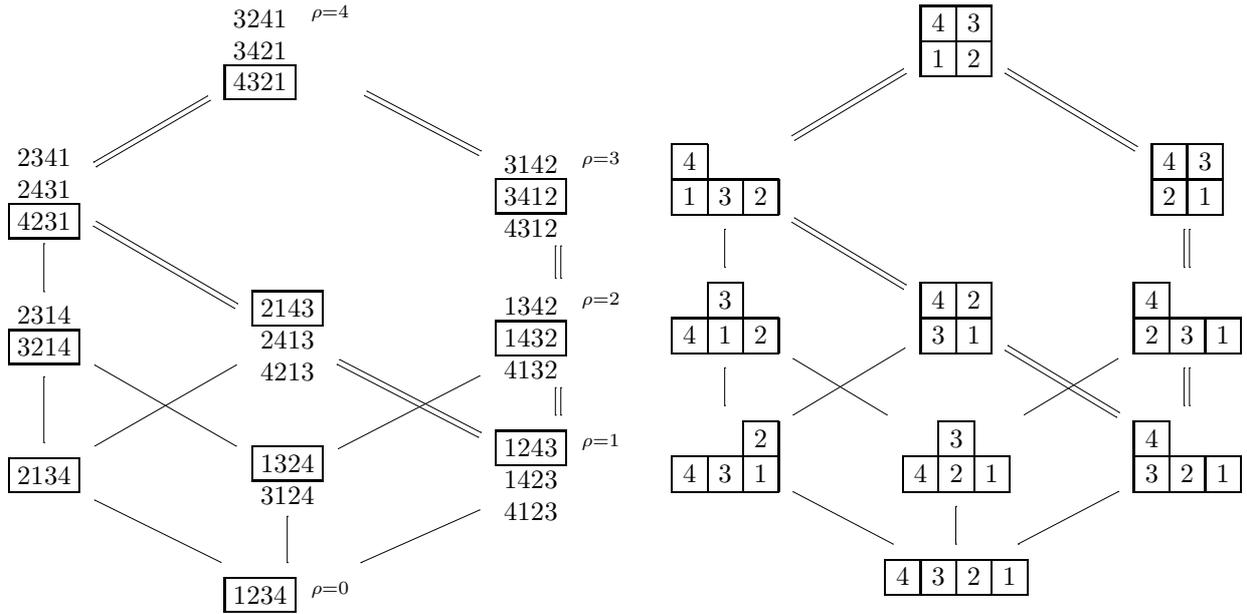

\caption{The graded weak order on Young-Fibonacci classes of size 4.}
\end{figure}

It is easy to check that each class appearing as a vertex of the
poset $(\youngFiboTab{4}, \younFiboTableauOrder)$ is an interval of
the weak order $(\mySn{4}, \permutohedronOrder)$, and this is a
general observation.

\begin{theorem}\label{theorem::classes::intervals}Let $t$ be a standard Young-Fibonacci tableau of size $n$,
then $\fiboClass{t}$ is an interval of the weak order $(\Sn,
\permutohedronOrder)$, more over $\fiboClass{t} = [\minCano{t},
\maxCano{t}]$.
\end{theorem}

To prove this statement, we will first relate $\fiboClass{t}$ with
linear extensions of a poset canonically associated to $t$, and then
we will prove that the set of linear extensions of this poset is an
interval of the weak order.

\begin{definition}\label{def::cano::poset}Let $t$ be a standard Young-Fibonacci tableau of size $n$, its
canonical poset $\canoposet{t}$ is the poset defined on the set
$\nn$ with the covering relations below.
\begin{enumerate}
  \item \label{cond::one::poset}the right-to-left reading of the bottom row of $t$ forms a chain in the
  poset ;
  \item each entry on top in a two-boxed column of $t$ is covered by
  the corresponding entry on bottom row.
\end{enumerate}
\end{definition}

\begin{note}
A permutation $\sigma$ is a \emph{toset} (totally ordered set) with
covering relations defined by $\sigma(i) \leq_\sigma \sigma(j)$
whenever $i < j$, that is to say $x \leq_\sigma y$ if $x$ appears to
the left of $y$ in $\sigma$. Let $\posetp$ be a poset and $\sigma$ a
permutation, $\sigma$ is said to be a linear extension of $\posetp$
if its relations preserve the relations in $\posetp$, that is to say
if $x \leq_{\posetp} y$ then $x \leq_{\sigma} y$. The set of linear
extensions of a poset $\posetp$ will be denoted
$\linearext{\posetp}$.
\end{note}

\begin{theorem}\label{theo::class::linear::extension}Let $t$ be a standard Young-Fibonacci tableau,
then $\fiboClass{t} = \linearext{\canoposet{t}}$.
\end{theorem}

\begin{proof}That any permutation $\sigma$ having $t$ as insertion
tableau is a linear extension of $\canoposet{t}$ is clear from
Definitions \ref{def::young::fibo::algo} and \ref{def::cano::poset}.
Conversely, if $\sigma$ is a linear extension of $\canoposet{t}$,
then $t$ is naturally built reading $\sigma$ from right to left
following the description given in Definition
\ref{def::young::fibo::algo}. At each new step the first letter one
reads is the maximal one (for $\leq_{\canoposet{t}}$) not yet read
in the chain described in rule (\ref{cond::one::poset}) of
Definition \ref{def::cano::poset}.
\end{proof}

\begin{theorem}\label{theo::extensions::canonical::poset}Let $t$ be a standard
$\youngFiboTableau$ of size $n$, then $\linearext{\canoposet{t}}$ is
the interval $[\minCano{t}, \maxCano{t}]$ in $(\Sn,
\permutohedronOrder)$.
\end{theorem}

To prove this statement we make use of the following well known
lemma.
\begin{lemma}\label{lemma::order::permuto}Let $\sigma$ and $\tau$ be
two permutations of $\Sn$, then the three properties below are
equivalent.
\begin{enumerate}
  \item \label{lemma::order::permuto::equiv::propOne} $\sigma \permutohedronOrder \tau$ ;
  \item \label{lemma::order::permuto::equiv::propTwo}
  $\ordSet{\tau} \subseteq \ordSet{\sigma}$ ;
  \item \label{lemma::order::permuto::equiv::propThree}
  $\invSet{\sigma} \subseteq \invSet{\tau}$.
\end{enumerate}
\end{lemma}

\def\myspace{\,\,}

\begin{proof} (of Theorem \ref{theo::extensions::canonical::poset})
It easily follows from the definition that $\canoposet{t}$ can be
partitioned into an antichain $A = (y_1, y_2, \cdots, y_\ell)$ and a
chain $C = (x_1 <_{\canoposet{t}} x_2 <_{\canoposet{t}} \cdots
<_{\canoposet{t}} x_k)$ such that for $i = 1 .. \ell$ there exists
$j(i) \leq k$ such that $y_i <_{\canoposet{t}} x_{j(i)}$, and
additionally for $i_1 < i_2$ one has $y_{i_1} < y_{i_2}$ and
$x_{j(i_1)} <_{\canoposet{t}} x_{j(i_2)}$. For illustrations, we use
the following example. \vspace{-0.25cm}
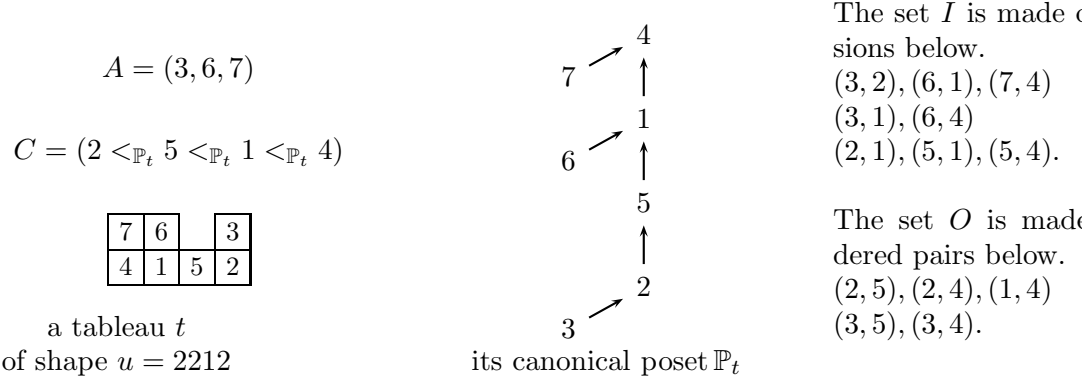
\begin{multicols}{3}
$$A = (3, 6, 7)$$
$$C = (2 <_{\canoposet{t}} 5 <_{\canoposet{t}} 1 <_{\canoposet{t}}
4)$$ {\small $$\ruban{7&6&\ &3\\4&1&5&2\\}$$}
\begin{tabular}{c}
a tableau $t$\\
of shape $u = 2212$
\end{tabular}

\columnbreak

\def\eUn{1} \def\eDeux{2} \def\eTrois{3} \def\eQuatre{4}
\def\eCinq{5} \def\eSix{6} \def\eSept{7}

\def\myposet{
\begin{psmatrix}[colsep=0.4,rowsep=0.1]
&& [name=e4]\eQuatre\\[0pt]
[name=e7]\eSept\\[0pt]
&& [name=e1]\eUn\\[0pt]
[name=e6]\eSix\\[0pt]
&& [name=e5]\eCinq\\[0pt]
\\[0pt]
&& [name=e2]\eDeux\\[0pt]
[name=e3]\eTrois
  \psset{nodesep=5pt,arrows=->}
  \ncline{e3}{e2} \ncline{e2}{e5} \ncline{e5}{e1} \ncline{e6}{e1} \ncline{e1}{e4} \ncline{e7}{e4}
\end{psmatrix}
}

$$\myposet$$
$$\mbox{its canonical poset}\, \canoposet{t}$$

\columnbreak

\noindent The set $I$ is made of the inversions below.\\$(3,2),
(6,1),
(7,4)$\\$(3,1), (6,4)$\\$(2,1), (5,1), (5,4)$.\\

\noindent The set $O$ is made of the ordered pairs below.\\$(2,5),
(2,4), (1,4)$\\$(3,5), (3,4)$.

\end{multicols}

\def\myspace{\,\,}
\noindent For $\sigma \in \linearext{\posetp}$, $\invSet{\sigma}$
includes at least the set $$I \,\, = \,\, \left\{
\begin{array}{l}
(y_i, x_{j(i)}), i = 1 .. \ell\\
(y_i, x_r) \myspace / \myspace x_{j(i)} > x_r \myspace
\mbox{and} \myspace x_{j(i)} <_{\canoposet{t}} x_r\\
(x_i, x_j) \myspace / \myspace x_i
> x_j \myspace \mbox{and} \myspace x_i <_{\canoposet{t}} x_j
\end{array} \right\}
$$ which is nothing but $\invSet{\minCano{t}}$ ;
so by [Lemma \ref{lemma::order::permuto} -
(\ref{lemma::order::permuto::equiv::propThree})], $\minCano{t}
\permutohedronOrder \sigma$. Moreover, $\ordSet{\sigma}$ includes at
least the set $$O \,\, = \,\, \{ \myspace (y_i, x_r) \myspace /
\myspace x_{j(i)} <_{\canoposet{t}} x_r \myspace \} \myspace \cup
\myspace \{ \myspace(x_i, x_j) \myspace / \myspace x_i < x_j
\myspace \mbox{and} \myspace x_i <_{\canoposet{t}} x_j \myspace \}$$
which is nothing but $\ordSet{\maxCano{t}}$ ; so by [Lemma
\ref{lemma::order::permuto} -
(\ref{lemma::order::permuto::equiv::propTwo})], $\sigma
\permutohedronOrder \maxCano{t}$ and hence $\sigma \in [\minCano{t},
\maxCano{t}]$. Conversely, for $\sigma \in [\minCano{t},
\maxCano{t}]$, applying Lemma \ref{lemma::order::permuto} to
$\minCano{t}$, $\sigma$ and $\maxCano{t}$ it appears that $\sigma$
has the inversions $y_i \leq_{\sigma} x_{j(i)}$ for $i = 1 .. \ell$,
and the relations $x_1 <_\sigma x_2 <_\sigma \cdots <_\sigma x_k$.
So $P(\sigma) = t$ and hence $\sigma \in \linearext{\canoposet{t}}$.
\end{proof}

\begin{proof} (of Theorem \ref{theorem::classes::intervals}) Follows
from Theorem \ref{theo::class::linear::extension} and Theorem
\ref{theo::extensions::canonical::poset}.\end{proof}

\newcommand{\rowcanonicaltab}[1]{rT_{#1}}
\newcommand{\columncanonicaltab}[1]{cT_{#1}}
\newcommand{\minrankoftab}[1]{\rho_{min}^{#1}}
\newcommand{\maxrankoftab}[1]{\rho_{max}^{#1}}

\begin{definition}Let $u$ be a snakeshape of size $n$, the row
canonical tableau $\rowcanonicaltab{u}$ is the one such that
\begin{enumerate}
  \item top cells of $\rowcanonicaltab{u}$ are labeled with entries $n$, $n-1$,
  $\cdots$ from left to right ;
  \item bottom cells in two-boxed columns are labeled with entries $1$, $2$,
  $\cdots$ from left to right.
  \end{enumerate}
The column canonical tableau $\columncanonicaltab{u}$ is built by
labeling the cells of $u$ from right to left and bottom to top.
\end{definition}

\begin{lemma}\label{lemma::min::max::rank}Let $u$ be a snakeshape of size $n$,
then $\columncanonicaltab{u}$ (resp. $\rowcanonicaltab{u}$) is the
unique tableau of shape $u$ having minimal rank $\minrankoftab{u}$
(resp. maximal rank $\maxrankoftab{u}$) in the poset
$(\youngFiboTab{n}, \younFiboTableauOrder)$. For any snakeshape $u$,
$\minrankoftab{u}$ is the number of double-boxed columns of $u$ and
$\maxrankoftab{u}$ is obtained as follows. Label each bottom cell
with the number of double-boxed columns on its left and do the same
but add $1$ for each top cell of double-boxed columns of $u$.
$\maxrankoftab{u}$ is the sum of labels obtained.
\end{lemma}

\newcommand\lmaxt[1]{\textcolor{blue}{#1}}
\newcommand\lmint[1]{\textcolor{red}{#1}}

\def\myu{2212}
\def\anyvertex{\bullet}
\def\tabOneOne{\hat{0}}
\def\tabTwoOne{\anyvertex} \def\tabTwoTwo{\anyvertex} \def\tabTwoThree{\anyvertex}
\def\tabThreeOne{\anyvertex} \def\tabThreeTwo{\anyvertex}
\def\tabThreeThree{\ruban{\lmint{7}&\lmint{5}&\ &\lmint{2}\\\lmint{6}&\lmint{4}&\lmint{3}&\lmint{1}\\}}
\def\tabFourOne{\anyvertex} \def\tabFourTwo{\anyvertex} \def\tabFourThree{\anyvertex}
\def\tabFiveOne{\anyvertex} \def\tabFiveTwo{\anyvertex} \def\tabFiveThree{\anyvertex}
\def\tabSixOne{\anyvertex}
\def\tabSixTwo{\ruban{\lmaxt{7}&\lmaxt{6}&\ &\lmaxt{4}\\\lmaxt{1}&\lmaxt{2}&\lmaxt{5}&\lmaxt{3}\\}}
\def\tabSixThree{\anyvertex}
\def\tabSevenOne{\hat{1}}
\def\arcinterm{\arx3}

{\small \centerline{
\newdimen\vcadre\vcadre=0.3cm 
\newdimen\hcadre\hcadre=0.3cm 
\setlength\unitlength{0.75mm}
$\xymatrix@R=0.4cm@C=0.8cm{%
*{\rho = 12} & & *{\tabSevenOne}\arx1[d]\arx1[dl]\arx1[dr]\\
*{\rho = 11} & *{\tabSixOne}\arcinterm[d] & *{\tabSixTwo}\arcinterm[d]\arcinterm[dr] & *{\tabSixThree}\arcinterm[d]\arcinterm[dl] %
& *{\mbox{row canonical tableau}} \\
& *{\tabFiveOne}\arcinterm[d]\arcinterm[dr] & *{\tabFiveTwo}\arcinterm[d]\arcinterm[dr] & *{\tabFiveThree}\arcinterm[d]\\
& *{\tabFourOne}\arcinterm[d]\arcinterm[drr] & *{\tabFourTwo}\arcinterm[dl] & *{\tabFourThree}\arcinterm[d]\arcinterm[dl]\\
*{\rho = 3} & *{\tabThreeOne}\arx1[d]\arx1[dr] & *{\tabThreeTwo}\arx1[d]\arx1[dl] & *{\tabThreeThree}\arx1[d]\arx1[dl] %
& *{\mbox{column canonical tableau}}\\
& *{\tabTwoOne}\arx1[dr] & *{\tabTwoTwo}\arx1[d] & *{\tabTwoThree}\arx1[dl]\\
*{\rho = 0} & & *{\tabOneOne}
 }$ }}

\vspace{-0.5cm}
\begin{figure}[h]
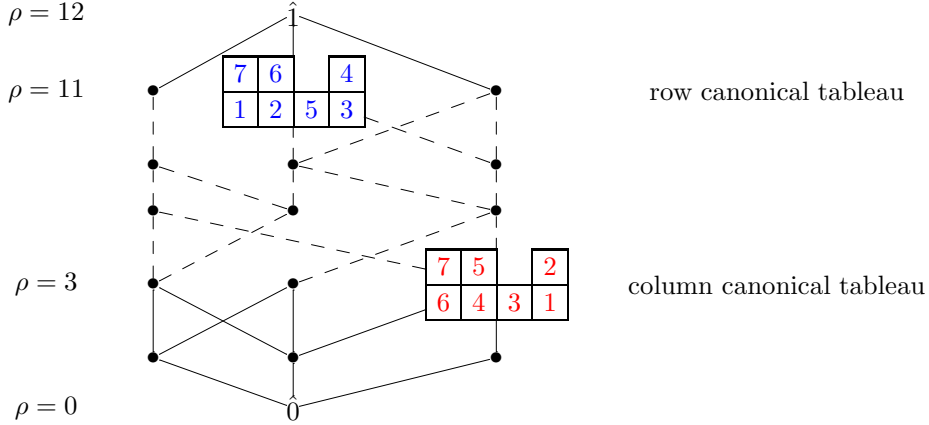

\caption{Row canonical and column canonical tableaux of shape
$2212$.}
\end{figure}


\begin{proof}(of Lemma \ref{lemma::min::max::rank}) Easily from the definitions.
\end{proof}

We will now relate $(\youngFiboTab{n}, \younFiboTableauOrder)$ to a
transition matrix in Okada's algebra associated to
$\youngFiboLattice$.

\section{A connection with Okada's algebra associated to the Young-Fibonacci lattice}\label{section::okada::algebra}

\newcommand{\okadakostkafibonumber}[2]{{\bf K}_{#1\!,\, #2}}

A Young-Fibonacci analogue of the ring of symmetric functions
\cite{macdo} was given and studied by S. Okada
\cite{okada_algebra_fibo}, with a Young-Fibonacci analogue of Kostka
numbers, appearing when expressing the analogue of a complete
symmetric function ${\bf h}_v$ in terms of the analogue of Schur
functions ${\bf s}_u$.
\begin{equation}
{\bf h}_v \,\, = \,\, \sum_u \okadakostkafibonumber{u}{v} \, {\bf
s}_u
\end{equation}

Young-Fibonacci analogue of Kostka numbers are generated by the
recurrence formulas below \cite{okada_algebra_fibo}, where
$\okadakostkafibonumber{a}{b}$ is defined for two snakeshapes of the
same weight and $\covers$ denotes the covering relation in
$\youngFiboLattice$.
\begin{equation}\label{recurrence::okada::fibo}
\left\{
\begin{array}{lclr}
\okadakostkafibonumber{1u}{1v} & = & \okadakostkafibonumber{u}{v} & \quad (r_1)\\
\okadakostkafibonumber{2u}{2v} & = & \okadakostkafibonumber{u}{v} & (r_2)\\
\okadakostkafibonumber{1u}{2v} & = & 0 & (r_3)\\
\okadakostkafibonumber{2u}{1v} & = & \sum_{w \covers u}
\okadakostkafibonumber{w}{v} & (r_4)
\end{array}\right.
\end{equation}

As it is stated below, the hook-length formula for binary trees
illustrated in Example \ref{example::hook} is an alternative formula
for computing $\okadakostkafibonumber{u}{1^n} =
\fiboTableauNumber{u}$ which is the dimension of a representation in
Okada's algebra.

\begin{proposition}Let $u$ be a snakeshape of size $n$, then
$\fiboTableauNumber{u}$ is the dimension of the module $V_u$
corresponding to $u$ in the $n^{th}$ homogenous component of Okada's
algebra associated to $\youngFiboLattice$.
\end{proposition}

\proof $dim(V_u)$ is the number of saturated chains from $\emptyset$
to $u$ in $\youngFiboLattice$, hence the result.\endproof

Here is a more general statement giving a combinatorial
interpretation of $\okadakostkafibonumber{u}{v}$ using
$(\youngFiboTab{n}, \younFiboTableauOrder)$.

\begin{theorem}\label{theo::poset::tableaux::kostka}Let $u$ and $v$ be two snakeshapes of size
$n$, and let $\hat{1}$ be the maximal tableau in $(\youngFiboTab{n},
\younFiboTableauOrder)$, then $\okadakostkafibonumber{u}{v}$ is the
number of tableaux $t$ of shape $u$ in the interval
$[\rowcanonicaltab{v}, \hat{1}]$.
\end{theorem}

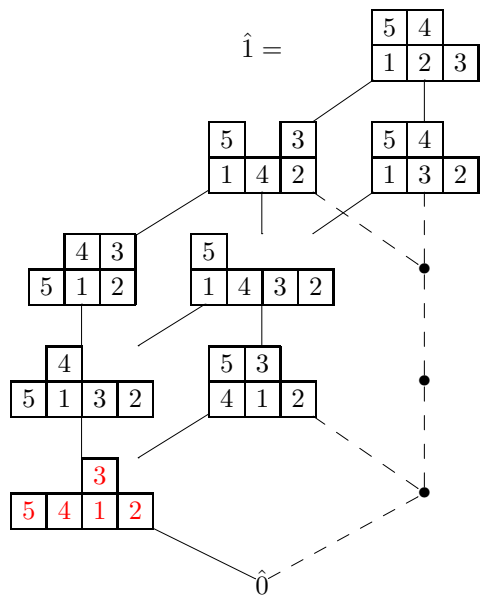
\begin{multicols}{2}

\newcommand\lgt[1]{\textcolor{blue}{#1}}
\newcommand\lbt[1]{\textcolor{red}{#1}}
\newcommand\canoelt[1]{\textcolor{red}{#1}}

\def\anyvertex{\bullet}
\def\tabOneOne{\hat{0}}
\def\tabTwoOne{\ruban{\ &\ &\canoelt{3}\\\canoelt{5}&\canoelt{4}&\canoelt{1}&\canoelt{2}\\}}
\def\tabTwoTwo{\anyvertex}
\def\tabThreeOne{\ruban{\ &4\\5&1&3&2\\}}
\def\tabThreeTwo{\ruban{5&3\\4&1&2\\}}
\def\tabThreeThree{\anyvertex}
\def\tabFourOne{\ruban{\ &4&3\\5&1&2\\}}
\def\tabFourTwo{\ruban{5\\1&4&3&2\\}}
\def\tabFourThree{\anyvertex}
\def\tabFiveOne{\ruban{5&\ &3\\1&4&2\\}}
\def\tabFiveTwo{\ruban{5&4\\1&3&2\\}}
\def\tabSixOne{\ruban{5&4\\1&2&3\\}}

{\small \centerline{
\newdimen\vcadre\vcadre=0.3cm 
\newdimen\hcadre\hcadre=0.3cm 
\setlength\unitlength{0.75mm}
$\xymatrix@R=0.55cm@C=0.5cm{%
& *{\hat{1} = } & *{\tabSixOne}\arx1[d]\arx1[dl] \\
& *{\tabFiveOne}\arx1[d]\arx1[dl]\arx3[dr] & *{\tabFiveTwo}\arx1[dl]\arx3[d]\\
*{\tabFourOne}\arx1[d] & *{\tabFourTwo}\arx1[d]\arx1[dl] & *{\tabFourThree}\arx3[d]\\
*{\tabThreeOne}\arx1[d] & *{\tabThreeTwo}\arx1[dl]\arx3[dr] & *{\tabThreeThree}\arx3[d] \\
*{\tabTwoOne}\arx1[dr] & & *{\tabTwoTwo}\arx3[dl] \\
& *{\tabOneOne}
 }$ }}

\columnbreak

\def\myzero{.}
\def\myv{1121}

\begin{example}In the matrix below, the number
$\okadakostkafibonumber{u}{\myv}$ counts the number of standard
Young-Fibonacci tableaux of shape $u$ in the interval
$[\rowcanonicaltab{\myv}, \hat{1}]$.\end{example}

\noindent {\small
\begin{tabular}{r|cccccccc}
& 221 & 212 & 2111 & 122 & 1211 & \canoelt{1121} & 1112 & $1^5$\\
\hline\\
221   & 1 & 1 & 2 & 1 & 2 & \canoelt{3} & 4 & 8\\
212   & \myzero & 1 & 1 & 1 & 1 & \canoelt{1} & 3 & 4\\
2111  & \myzero & \myzero & 1 & \myzero & 1 & \canoelt{1} & 1 & 4\\
122   & \myzero & \myzero & \myzero & 1 & 1 & \canoelt{1} & 2 & 3\\
1211  & \myzero & \myzero & \myzero & \myzero & 1 & \canoelt{1} & 1 & 3\\
1121  & \myzero & \myzero & \myzero & \myzero & \myzero & \canoelt{1} & 1 & 2\\
1112  & \myzero & \myzero & \myzero & \myzero & \myzero & \myzero & 1 & 1\\
$1^5$ & \myzero & \myzero & \myzero & \myzero & \myzero & \myzero &
\myzero & 1
\end{tabular}
}

Iterating this for each snakeshape $v$ of size $n$, one builds the
transition matrix for expressing the analogue of complete symmetric
function ${\bf h}_v$ in terms of the analogue of Schur functions
${\bf s}_u$.
\end{multicols}

\vspace{-0.5cm}
\begin{figure}[h]\label{fig::kostka::fibo::numbers::matrix}
\caption{$(\youngFiboTab{5}, \younFiboTableauOrder)$ and Okada's
analogue of Kostka matrix for $n=5$.}
\end{figure}

\begin{proof}(of Theorem \ref{theo::poset::tableaux::kostka}) A proof consists in showing that for any couple $(a, b)$
of snakeshapes appearing in the left hand side of equation
(\ref{recurrence::okada::fibo}), there is a one-to-one
correspondence between tableaux satisfying the conditions of the
theorem for $(a,b)$ and those satisfying the conditions of the
theorem for the couples of snakeshapes in the corresponding right
hand side of the relation. For $(r_1)$, given a tableau $t$ of shape
$u$ such that $\rowcanonicaltab{v} \younFiboTableauOrder t$, $t$ is
mapped onto the tableau $t'$ of shape $1u$ obtained from $t$ by
attaching a cell labeled $n+1$ to its left, and
$\rowcanonicaltab{1v} \younFiboTableauOrder t'$. For $(r_2)$, one
attaches a two-boxed column to the left of $t$, with $1$ as bottom
entry and $n+2$ as top entry, in addition one standardizes $t$ by
increasing all its entries. Then $t'$ is of shape $2u$ and
$\rowcanonicaltab{2v} \younFiboTableauOrder t'$. For $(r_3)$ it
easily follows from the definition of the operation of shifting an
entry in a tableau that there is no tableaux $t_1$ and $t_2$ of
shape $1u$ and $2v$ respectively, such that $t_2
\younFiboTableauOrder t_1$. For $(r_4)$, let $t$ be a tableau of
shape $2u$ such that $\rowcanonicaltab{1v} \younFiboTableauOrder t$,
then $t$ is mapped onto the tableau $t' =
\evacuatedTableauFromLetter{t}{n}$, that is the tableau obtained
from $t$ by evacuating its maximal letter (the evacuation process
originally due to Killpatrick \cite{evacuation_fibonacci} is
described in Section \ref{section::discussion::evacuation}). Indeed,
let $w$ be the shape of $t'$, then $w \covers u$ and
$\rowcanonicaltab{v} \younFiboTableauOrder t'$.
\end{proof}
\section{Kostka numbers, the Littlewood Richardson rule, and four posets on Young tableaux}
\label{section::kostka::lr::young::tableaux}

\newcommand\restrictedTabPartit[3]{\lambda(#1_{/#2, #3})}

The poset $(\youngFiboTab{n}, \younFiboTableauOrder)$ of
Young-Fibonacci tableaux we defined in Section
\ref{section::order::fibo::tableaux} is an analogue of one among
four partial orders on the set $\youngTab{n}$ of standard Young
tableaux of size $n$ \cite{taskin}. The weak order $(\youngTab{n},
\younTableauWeakOrder)$ is defined as in Theorem
\ref{theo::weakest::order} with $P(\sigma)$ denoting the Schensted
insertion tableau of $\sigma$. Let $\lambda$ and $\mu$ be two
partitions of lengths $\ell(\lambda)$ and $\ell(\mu)$, $\lambda$ is
said greater than $\mu$ in the dominance order and one writes
$\lambda
\partitDominanceOrder \mu$ if for each $1 \leq i \leq min(\ell(\lambda), \ell(\mu))$,
the inequality $\lambda_1 + \lambda_2 + \cdots + \lambda_i \geq
\mu_1 + \mu_2 + \cdots + \mu_i$ holds. Let $t$ be a standard Young
tableau of size $n$, and $1 \leq i \leq j \leq n$. We denote
$\restrictedTabPartit{t}{i}{j}$ the shape of the tableau obtained
from $t$ by first restricting $t$ to the segment $[i,j]$, then
lowering all entries by $i-1$, and finally sliding the skew tableau
obtained into normal shape by jeu-de-taquin. The chain order
$\younTableauChainOrder$ on standard Young tableaux is defined as
follows.

\begin{definition}\label{def::order::chain::young::tableau}\cite{taskin} Let $t$ and $t'$ be two standard
Young tableaux of size $n$, then $t \younTableauChainOrder t'$ if
and only if for each $1 \leq i \leq j \leq n$,
$\restrictedTabPartit{t}{i}{j} \partitDominanceOrder
\restrictedTabPartit{t'}{i}{j}$.
\end{definition}

The reader interested may refer to \cite{taskin} for the definition
of the two other orders, as well as for the properties of those
posets. The four posets happen to coincide up to rank $n = 5$.

\def\matrixColSep{0}
\def\matrixRowSep{-0.1}

\begin{center}
$\xymatrix@R=1.3cm@C=0.6cm{%
 &  &  &  &  &  & *{\ytVingtSix}& \\
 &  & *{\ytVingtDeux}\arx1[rrrru]&  &  & *{\ytVingtTrois}\arx1[ru]&  & *{\ytVingtQuatre}\arx1[lu]&  &  &
*{\ytVingtCinq}\arx1[llllu]& \\
 & *{\ytDixSept}\arx1[ru]\arx1[rrrru]&  &  & *{\ytDixHuit}\arx1[ru]&  &
*{\ytDixNeuf}\arx1[llllu]\arx1[rrrru]&  &  & *{\ytVingt}\arx1[llu]&
&
*{\ytVingtEtUn}\arx1[llllu]\arx1[lu]& \\
*{\ytOnze}\arx1[ru]&  &  &  & *{\ytDouze}\arx1[lluu]\arx1[rrrrru]& &
*{\ytTreize}\arx1[llu]\arx1[rrru]&  & *{\ytQuatorze}\arx1[llu]&
*{\ytQuinze}\arx1[lllllu]\arx1[ruu]&  &  &
*{\ytSeize}\arx1[lu]& \\
 & *{\ytSix}\arx1[lu]\arx1[rrruu]&  &  & *{\ytSept}\arx1[llluu]\arx1[rrrru]\arx1[rrrrru]&  &
*{\ytHuit}\arx1[u]&  &  &
*{\ytNeuf}\arx1[lllllu]\arx1[lu]\arx1[rruu]& &
*{\ytDix}\arx1[lluu]\arx1[ru]& \\
 &  & *{\ytDeux}\arx1[lu]\arx1[rruu]\arx1[rrrru]&  &  & *{\ytTrois}\arx1[llllu]\arx1[lu]&  &
*{\ytQuatre}\arx1[rru]\arx1[rrrru]&  &  & *{\ytCinq}\arx1[llllu]\arx1[luu]\arx1[ru]& \\
 &  &  &  &  &  & *{\ytUn}\arx1[llllu]\arx1[lu]\arx1[ru]\arx1[rrrru]& \\
}$\end{center}

\vspace{-0.75cm}
\begin{figure}[h]
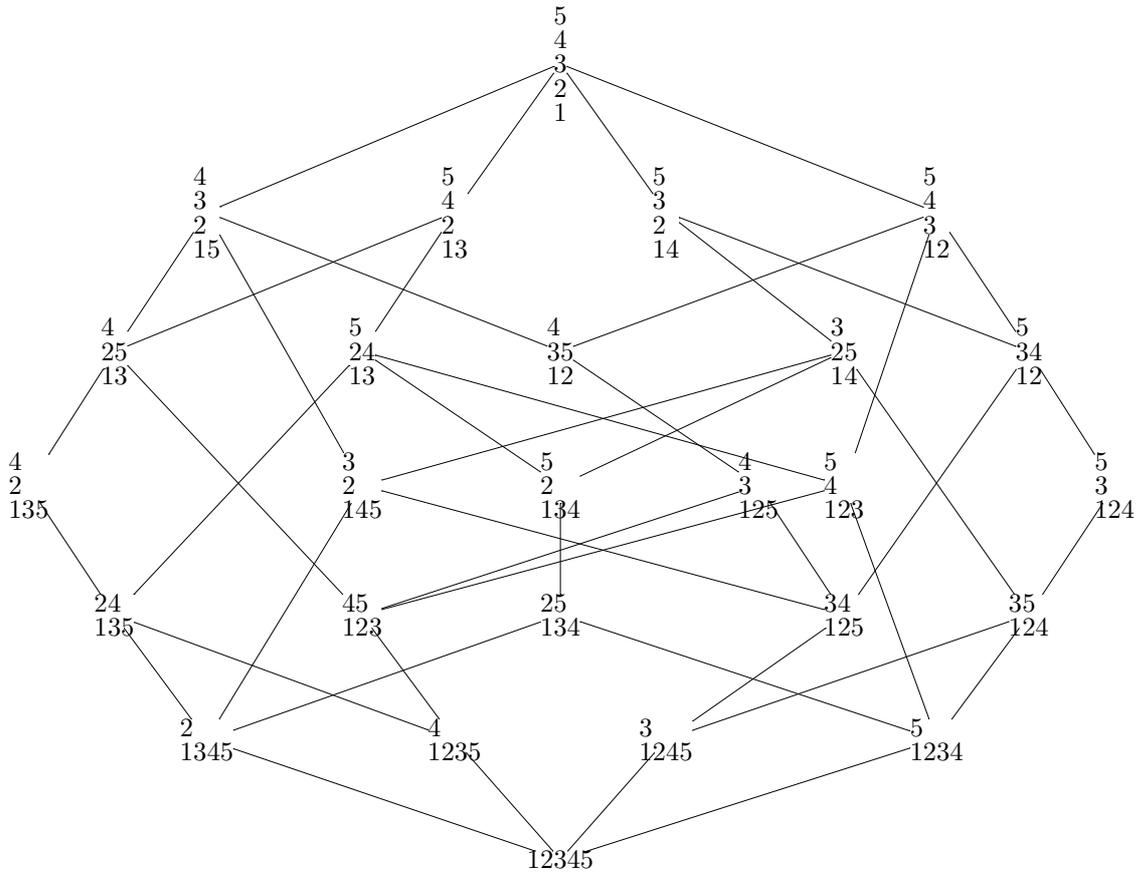

\caption{Partial order on Young tableaux of size 5.}
\end{figure}

Below is a Young tableaux analogue of Theorem
\ref{theo::poset::tableaux::kostka}.

\def\nsclt{nsclt(n)}
\def\nsrtmu{nscrt(\mu)}

\begin{theorem}\label{theo::poset::young::tableaux::kostka}Let $\lambda$, $\mu$ be two partitions of size $n$,
let $\rowcanonicaltab{\mu}$ be the row canonical standard Young
tableau of shape $\mu$, that is to say $\rowcanonicaltab{\mu}$ has
shape $\mu$ and is increasingly filled from let to right and bottom
to top. And let $\hat{0}$ be the minimal tableau in the poset of
standard Young tableaux of size $n$. Then
$\kostkanumber{\lambda}{\mu}$ is the number of standard Young
tableaux of shape $\lambda$ in the interval $[\hat{0},
\rowcanonicaltab{\mu}]$, for any one of the posets studied in
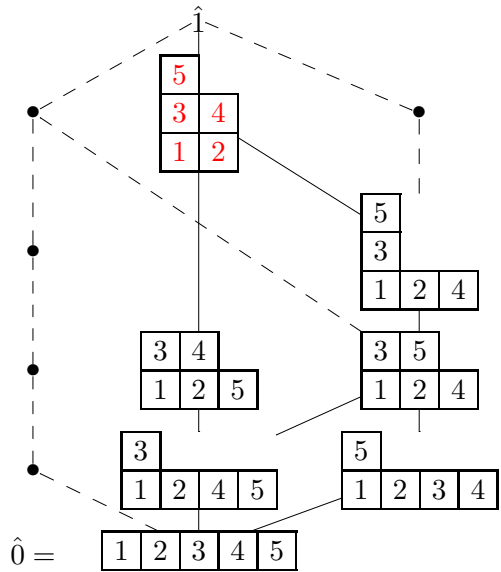
\cite{taskin}.
\end{theorem}

\begin{multicols}{2}

\newcommand\canoelt[1]{\textcolor{red}{#1}}

\def\anyvertex{\bullet}
\def\tabOneOne{\ruban{1&2&3&4&5\\}}
\def\tabTwoOne{\anyvertex}
\def\tabTwoTwo{\ruban{3\\1&2&4&5\\}}
\def\tabTwoThree{\ruban{5\\1&2&3&4\\}}
\def\tabThreeOne{\anyvertex}
\def\tabThreeTwo{\ruban{3&4\\1&2&5\\}}
\def\tabThreeThree{\ruban{3&5\\1&2&4\\}}
\def\tabFourOne{\anyvertex}
\def\tabFourTwo{\anyvertex}
\def\tabFourThree{\ruban{5\\3\\1&2&4\\}}
\def\tabFiveOne{\anyvertex}
\def\tabFiveTwo{\ruban{\canoelt{5}\\\canoelt{3}&\canoelt{4}\\\canoelt{1}&\canoelt{2}\\}}
\def\tabFiveThree{\anyvertex}

\def\tabSixOne{\hat{1}}

\centerline{
\newdimen\vcadre\vcadre=0.3cm 
\newdimen\hcadre\hcadre=0.3cm 
\setlength\unitlength{0.75mm}
$\xymatrix@R=0.3cm@C=0.6cm{%
& *{\tabSixOne}\arx3[d]\arx3[dl]\arx3[dr]\\
*{\tabFiveOne} & *{\tabFiveTwo} & *{\tabFiveThree}\\
*{\tabFourOne}\arx3[u] &  & *{\tabFourThree}\arx1[ul]\arx3[u]\\
*{\tabThreeOne}\arx3[u] & *{\tabThreeTwo}\arx1[uu] & *{\tabThreeThree}\arx1[u]\arx3[uull]\\
*{\tabTwoOne}\arx3[u] & *{\tabTwoTwo}\arx1[u]\arx1[ur] & *{\tabTwoThree}\arx1[u]\\
*{\hat{0} =} & *{\tabOneOne}\arx1[u]\arx1[ur]\arx3[ul]
 }$ }

\columnbreak

\def\myzero{.}

\begin{example}In the matrix below, the number
$\kostkanumber{\lambda}{221}$ counts the number of standard Young
tableaux of shape $\lambda$ in the interval $[\hat{0},
\rowcanonicaltab{221}]$.\end{example}

\begin{center}\begin{tabular}{rr|ccccccc}
& $\mu$ = & 5 & 41 & 32 & 311 & \canoelt{221} & 2111 & 11111\\
\hline\\
$\lambda$ = & 5   & 1 & 1 & 1 & 1 & \canoelt{1} & 1 & 1\\
& 41   & \myzero & 1 & 1 & 2 & \canoelt{2} & 3 & 4\\
& 32  & \myzero & \myzero & 1 & 1 & \canoelt{2} & 3 & 5\\
& 311   & \myzero & \myzero & \myzero & 1 & \canoelt{1} & 3 & 6\\
& 221  & \myzero & \myzero & \myzero & \myzero & \canoelt{1} & 2 & 5\\
& 2111  & \myzero & \myzero & \myzero & \myzero & \myzero & 1 & 4\\
& 11111 & \myzero & \myzero & \myzero & \myzero & \myzero & \myzero
& 1
\end{tabular}\end{center}

Iterating this for each partition $\mu$ of size $n$, one builds the
transition matrix for expressing the complete symmetric function
$h_\mu$ in terms of the Schur functions $s_\lambda$.

\end{multicols}

\vspace{-0.5cm}
\begin{figure}[h]\label{fig::poset::young::tableau::kostka}
\caption{Poset of Young tableaux and Kostka matrix for $n=5$.}
\end{figure}

\begin{proof}(of Theorem \ref{theo::poset::young::tableaux::kostka})
For a given partition $\mu$, let $\nsrtmu$ be the row canonical
semi-standard Young tableau of shape $\mu$, that is the tableau
filled with 1's on its first line, 2's on its second line and so on.
Let $\nsclt$ be the semi-standard Young tableau of shape $n$ and
having $\mu_i$ entries $i$ for $i = 1..\ell(\mu)$. Consider the
extension of Definition \ref{def::order::chain::young::tableau} to
the set $Tab(\mu)$ of semi-standard Young tableaux having $\mu_i$
entries $i$ for $i = 1..\ell(\mu)$. Then for each $t \in Tab(\mu)$,
one has $\nsclt \younTableauChainOrder t \younTableauChainOrder
\nsrtmu$. There is a canonical bijection mapping $(Tab(\mu),
\younTableauChainOrder)$ onto $([\hat{0}, \rowcanonicaltab{\mu}],
\younTableauChainOrder)$ and this map is order preserving. So
Theorem \ref{theo::poset::young::tableaux::kostka} holds for the
partial order $\younTableauChainOrder$. From (\cite{taskin}, Theorem
1.1) and the remark that $[\hat{0}, \rowcanonicaltab{\mu}] =
\rowcanonicaltab{\mu_1} * \rowcanonicaltab{\mu_2} * \cdots
* \rowcanonicaltab{\mu_{\ell(\mu)}}$, it follows that the set
of tableaux in $[\hat{0}, \rowcanonicaltab{\mu}]$ does not depend on
the choice of the partial order.
\end{proof}

\section*{Concluding remarks and perspectives}

There are quite many similarities between the Robinson-Schensted
algorithm and the Young-Fibonacci insertion algorithm. As well as
between the combinatorics of Young tableaux and the combinatorics of
Young-Fibonacci tableaux. One of the questions we have not explored
in this paper is the one of the existence of an algebra of
Young-Fibonacci tableaux, which would be an analogue of the
Poirier-Reutenauer Hopf algebra of Young tableaux
\cite{algebre_tableaux}. Such an algebra would certainly help in
giving a combinatorial description (in terms of tableaux) of the
product of Schur functions in Okada's algebra associated to the
Young-Fibonacci lattice. We are currently looking for a suitable
definition of this algebra.

\section*{Acknowledgements}
The author is grateful  to  F. Hivert for helpful comments and
suggestions throughout this work.

\bibliographystyle{amsalpha}

\begin{thebibliography}{A}

\def\sameauthor{---}

\bibitem{schensted_subsequence}{\sl C. Schensted}, Longest increasing and decreasing subsequences.
{\it Canad. J. Math.}, vol. 13, 1961, pp. 179-191.

\bibitem{knuth_rsk} {\sl D. E. Knuth},
Permutations, matrices and generalized Young tableaux, {\it Pacific
J. Math.} 34(1970), 709--727.

\bibitem{art_of_comp_programming}{\sl \sameauthor},
The art of computer programming, vol.3: {\it Searching and sorting}
(Addison-Wesley, 1973).

\bibitem{ncsf4}{\sl D. Krob} and {\sl J.-Y. Thibon},
Noncommutative symmetric function IV: Quantum linear groups and
Hecke algebras at q=0, {\it J. Alg. Comb}. 6 (1997), 339-376.

\bibitem{f_hivert_pbt}{\sl F. Hivert}, {\sl J. C. Novelli}, and
{\sl J.-Y. Thibon}, The Algebra of Binary Search Trees, {\it Theo.
Comp. Science} 339(2005), 129-165.

\bibitem{macdo}{\sl I. G. MacDonald}, Symmetric functions and Hall Polynomials, {\it 2nd ed, Clarendon
Press, Oxford Sce Publications}, 139(1995).

\bibitem{nzeutchap}{\sl J. Nzeutchap},
On the Young-Fibonacci Insertion Algorithm, to appear in the
Proceedings of {\it FPSAC'07}.

\bibitem{evacuation_fibonacci}{\sl K. Killpatrick},
Evacuation and a Geometric Consturction for Fibonacci Tableaux, {\it
J. Comb. Th, Ser A}, 110 (2005), 337-351.

\bibitem{taskin}{\sl M. Taskin},
Properties of four partial orders on standard Young tableaux, {\it
J. Comb. Theory, Ser A}, 113(2006), 1092-1119.

\bibitem{differential_posets}{\sl R. P. Stanley},
Differential Posets, {\it J. Amer. Math. Soc.} 1 (1998), 919-961.

\bibitem{fibolattice_lattice}{\sl \sameauthor},
The Fibonacci lattice, {\it Fibonacci Quarterly} 13(1998), 215-232.

\bibitem{algebre_tableaux}{\sl S. Poirier} and {\sl C. Reutenauer},
Alg\`ebre de Hopf des tableaux, {\it Ann. Sci. Math. Q\'ebec} 19
(1995), 79-90.

\bibitem{graph_duality}{\sl S. V. Fomin},
Duality of Graded Graphs, {\it J. Alg. Comb.} 3(1994), 357-404.


\bibitem{generalized_rsk}{\sl \sameauthor}, Generalized Robinson-Schensted-Knuth correspondence,
{\it Zapiski Nauchn. Sem. LOMI.} 155 (1986), 156-175.

\bibitem{schensted_for_dual_graphs}{\sl \sameauthor},
Schensted Algorithms for Dual Graded Graphs, {\it J. Alg. Comb.}
4(1995), 5-45.

\bibitem{okada_algebra_fibo}{\sl S. Okada},
Algebras associated to the Young-Fibonacci lattice, {\it Trans AMS}
346(1994), 549-568.

\bibitem{roby_thesis}{\sl T. Roby},
Applications and extensions of Fomin's generalization of the
Robinson-Schensted correspondence to differential posets, {\it Ph.D.
thesis, MIT}, 1991.

\end{thebibliography}

\end{document}